\documentclass[oneside,leqno,11pt]{article}
\usepackage{amssymb,amsmath,latexsym}

\newcommand{\s}{\smallskip}
\newcommand{\m}{\medskip}
\renewcommand{\b}{\bigskip}
\newcommand{\diag}{\operatorname{diag}}
 
\def\ga{\mathfrak{a}}
\def\gb{\mathfrak{b}}

\def\gg{\mathfrak{g}}
\def\gh{\mathfrak{h}}

\def\gj{\mathfrak{j}}
\def\gk{\mathfrak{k}}
\def\gl{\mathfrak{l}}
\def\gm{\mathfrak{m}}
\def\gn{\mathfrak{n}}

\def\gp{\mathfrak{p}}
\def\gq{\mathfrak{q}}
\def\gr{\mathfrak{r}}
\def\gs{\mathfrak{s}}
\def\gt{\mathfrak{t}}
\def\gu{\mathfrak{u}}

\def\gz{\mathfrak{z}}

\def\C{\mathbb{C}}

\def\E{\mathbb{E}}
\def\F{\mathbb{F}}

\def\H{\mathbb{H}}

\def\R{\mathbb{R}}

\def\Z{\mathbb{Z}}

\def\cB{\mathcal{B}}

\def\cD{\mathcal{D}}

\def\cF{\mathcal{F}}

\def\cH{\mathcal{H}}

\def\cK{\mathcal{K}}

\def\cO{\mathcal{O}}

\def\cU{\mathcal{U}}
\def\cV{\mathcal{V}}

\renewcommand{\thesection}{\arabic{section}}
\renewcommand{\theequation}{\thesection.\arabic{equation}}

\newtheorem{theorem}[equation]{Theorem}

\newtheorem{lemma}[equation]{Lemma}

\newtheorem{corollary}[equation]{Corollary}

\newtheorem{proposition}[equation]{Proposition}

\newtheorem{definition}[equation]{Definition}

\newtheorem{example}[equation]{Example}

\newtheorem{remark}[equation]{Remark}

\title{Principal Series Representations of Direct Limit Groups}

\date{9 February 2004}

\author{Joseph A. Wolf\,\footnote{
Research partially supported by NSF Grant DMS 99-88643.
\endgraf
\endgraf
}
}

\begin{document}

\maketitle

\abstract{\begin{quote} \footnotesize 
We combine the geometric realization of principal series representations
of \cite{W1} with the Bott--Borel--Weil Theorem for direct limits of
compact groups found in \cite{NRW3}, obtaining limits of principal series 
representations for direct limits of real reductive Lie groups.  We
introduce the notion of weakly parabolic direct limits and relate it to
the conditions that the limit representations are norm--preserving 
representations on a Banach space or unitary representations on a Hilbert
space.  We specialize the results to diagonal embedding direct limit groups.
Finally we discuss the possibilities of extending the results to limits of
tempered series other than the principal series.
\end{quote} }
\vskip 1 cm

\section{Introduction} \label{sec1}
\setcounter{equation}{0}

Harmonic analysis on a real reductive Lie group $G$ depends on several series
of representations, one for each conjugacy class of Cartan subgroups of $G$.  
See \cite{H0}, \cite{H1}, \cite{H2} and \cite{H3} for the case where
$G$ is Harish--Chandra class, \cite{W1}, \cite{HW1} and \cite{HW2} for
the general case.  The simplest of these series is the {\em principal
series}.  It consists of representations constructed from representations
of compact Lie groups, characters on real vector groups, and the induced
representation construction.  The other series are somewhat more delicate,
replacing \' Elie Cartan's theory of representations of compact Lie groups 
by Harish--Chandra's theory of discrete series representations of real 
reductive Lie groups.
\m

This paper is the first step in a program to extend the construction,
analysis and geometry of those series of representations from the
finite--dimensional setting to a nontrivial
but well behaved family of infinite--dimensional Lie groups, the
direct limits of real reductive Lie groups.  Here we consider the case
of the principal series.  The case of the discrete series, and then the 
general case, will be considered separately in \cite{W5} and \cite{W6}.
\m

The classical Bott--Borel--Weil Theorem \cite{Bo} realizes representations
of compact Lie groups as cohomology spaces of holomorphic vector bundles
over complex flag manifolds.  It since has been extended to direct limits
of compact Lie groups and direct limits of complex Lie groups, both in the 
analytic category \cite{NRW3} and in the algebraic category \cite{DPW}.
With some technical adjustment, the results of \cite{NRW3} replace
Cartan's theory of representations of compact Lie groups for construction
of direct limit principal series representations.  There are, however, a
number of technical points, some of them delicate, that have to be addressed
and we mention them as we describe the contents of this paper.
\m

Section \ref{sec2} recalls our class of finite--dimensional real reductive
Lie groups and the standard construction of their not--necessarily--unitary
principal series representations.  Section \ref{sec3} recalls the
geometric realization of those representations on partially holomorphic
cohomologies of vector bundles over closed orbits in complex flag manifolds.
In Section \ref{sec4} we discuss alignment questions for minimal
parabolic subgroups.  The alignment is needed in order to define limit
principal series representations of our direct limit groups.  In effect,
this is the first technical issue, and it addresses the question of
whether $G = \varinjlim G_i$ can have a meaningful
direct limit of principal series representations.  For that we need
the connecting maps $\phi_{j,i}:G_i \to G_j$ of the
direct system to respect the ingredients of the principal series recipe.
Initially that must be done for the components $M_i$,
$A_i$ and $N_i$ of minimal parabolic subgroups $P_i =
M_i A_i N_i \subset G_i$\,.  That alignment on components
is not quite automatic, but it holds (possibly after passing to a cofinal
subsystem --- which yields the same limit group) for the most interesting
cases, the diagonal embedding direct limit groups of Section \ref{sec9}.  See
Proposition \ref{diag-lim-exists}.  Next, it must be done on the level
of representations of the $M_i$\,.  That, of course, is automatic for
spherical principal series representations, but more generally
we use an appropriate extension
of Cartan's highest weight theory.  Thus we obtain representations of $G$ that
are direct limits of principal series representations of the $G_i$\,.
\m
 
The second issue is to construct good geometric realizations of these
``principal series'' representations of the limit groups $G$.  This is the
heart of the paper.
The method of \cite{W2}, illustrated in \cite[Section 1]{W2}, gives natural
partially holomorphic realizations of principal series representations
$\pi_i$ of $G_i$\,.  That involves a certain extension of the
classical Bott--Borel--Weil Theorem \cite{Bo} which we need for the
groups $M_i$\,.  In order to pass to the limit, we construct and study
the appropriate limit flag manifolds, limit of closed orbits, limits of
holomorphic arc components, and limit sheaves, in Sections \ref{sec5} and
\ref{sec6}.  This is done in such a way that the limit Bott--Borel--Weil
Theorem of \cite{NRW3} applies over the  holomorphic arc components of the
closed orbits.  That defines the geometric setting for the representations
in question.  In order to see that the cohomology of the limit sheaf
is the limit of the cohomologies, we prove a Mittag--Leffler condition
at the end of Section \ref{sec6}.  Thus we have the possibility of obtaining
good geometric realizations of limit principal series representations of $G$
directly on cohomology spaces.
\m

We actually construct the geometric realizations in Section \ref{sec7}.
Theorem \ref{limit-bw} is the $0$--cohomology result in the style of the
Borel--Weil Theorem, and Theorem \ref{limit-bbw} is the higher cohomology
result in the style of the Bott--Borel--Weil Theorem.  For the latter
it is essential to have the cohomologies all occur in the same degree.
That is the third technical issue, and we reduce it to the same question
for $M = \varinjlim M_i$\,, where it was settled in \cite{NRW3}.
\m
 
The fourth issue is whether these principal series representations of
$G$ are norm--preserving Banach space representations, or even unitary
representations, of $G$.  That is settled in Theorem \ref{lp-limit}.  There
the key idea is that of weakly parabolic direct system.
\m
 
It is very important to have a large number of interesting examples.
For that we consider diagonal embedding direct limits of classical real
simple Lie groups.  We examine their behavior relative to the various
general notions studied earlier and see that our constructions work very
well for these interesting direct limit groups.  This is done in
Section \ref{sec9}.  These diagonal
embedding direct limits have been studied extensively in the context of
locally finite Lie algebras. That is a rapidly developing area; see
\cite{LN} and the references there.  A
locally finite Lie algebra of countable dimension can be represented as
a direct limit $\varinjlim \{\gg_m,d\phi_{n,m}\}_{m,n \in \Z^+}$ of finite
dimensional Lie algebras, and the diagonal embedding direct limits are
essentially just those where the group level maps $\phi_{n,m}$ are
polynomials of degree $1$.
\m
 
Finally, in Section \ref{sec10} we
discuss the place of the principal series in our program for constructing
limit representations corresponding to all tempered series, and indicate
some of the problems to be settled in \cite{W5} and \cite{W6}.
\m

The notions of parabolic and weakly parabolic direct systems developed
from a conversation with Andrew Sinton.

\section{Principal Series for General Reductive Groups} \label{sec2}
\setcounter{equation}{0}

Let $G$ be a reductive real Lie group.  In other words, its Lie
algebra $\gg$ is reductive in the sense that it is the direct sum
of a semisimple Lie algebra $\gg' = [\gg , \gg]$ and an abelian
idea $\gz$ which is the center of $\gg$\,.  As usual, $\gg_{_\C}$
denotes the complexification of $\gg$\,, 
so $\gg_{_\C} = \gg_{_\C}' \oplus \gz_{_\C}$
direct sum of the respective complexifications of $\gg'$ and
$\gz$\,.  We assume that $G$ satisfies the conditions of \cite{W1}:
\begin{equation} \label{group-class1}
\text{if } g \in G \text{ then } Ad(g) \text{ is an inner
automorphism of } \gg_{_\C} \text{\,, and } 
\end{equation}
\begin{equation} \label{group-class2}
\begin{aligned}
&G \text{ has a closed normal abelian subgroup } Z \text{ such that }\\
&\phantom{XXX} Z \text{ centralizes the identity component } 
	G^0 \text{ of } G \ , \\
&\phantom{XXX} ZG^0 \text{ has finite index in } G \ \text{, and} \\
&\phantom{XXX} Z \cap Z_{G^0} \text{ is co-compact in the center } 
	Z_{G^0} \text{ of } G^0 \ .
\end{aligned}
\end{equation}
These are the conditions, inherited by Levi components of cuspidal
parabolic subgroups of $G$\,, that lead to a nice Plancherel formula.
See \cite{W1}, \cite{HW1}, and \cite{HW2}.  The famous {\em Harish--Chandra
class} is the case where the semisimple component $(G^0)' := [G^0, G^0]$ of
$G^0$ has finite center and the component group $G/G^0$ is finite.
\m

Condition (\ref{group-class1}) says that the standard tempered representation
construction yields representations that have an infinitesimal character.
It can be formulated: $Ad(G) \subset Int(\gg_{_\C})$.  
\s

Note that the kernel of
$Ad : G \to Ad(G)$ is the centralizer $Z_G(G^0)$ of the identity
component and that the image $Ad(G)$ is a closed subgroup of the 
complex semisimple group $Int(\gg_{_\C})$ with only finitely many topological 
components.  Thus $Ad(G)$ has maximal compact subgroups, as usual for
semisimple linear groups, and every maximal compact subgroup of $Ad(G)$ 
is of the form $K/Z_G(G^0)$ for some closed subgroup 
$K \subset G$\,.  
\s

Given a maximal compact subgroup $K/Z_G(G^0)$ of $Ad(G)$, 
it is known \cite[Lemma 4.1.1]{W1} 
that $K$ is the fixed point set of a unique involutive automorphism
$\theta$ of $G$\,.  These automorphisms $\theta$ are called {\em Cartan
involutions} of $G$\,, and they are lifts of the Cartan involutions
of the linear group $Ad(G)$.  The groups $K$ are the {\em maximal
compactly embedded} subgroups of $G$.
\s

One also knows \cite[Lemma 4.1.2]{W1} that $K \cap G^0$ is the identity
component $K^0$ of $K$\,, that $K$ meets every topological component
of $G$, that any two Cartan involutions of $G$ are conjugate by an
element of $Ad(G^0)$\,, and that every Cartan subgroup of $G$ is stable
under some Cartan involution.  Here we use the usual definition: Cartan
subgroup of $G$ means the centralizer of a Cartan subalgebra of $\gg$.
\m

Fix a Cartan involution $\theta$ of $G$ and the corresponding
maximal compactly embedded subgroup $K = G^\theta$ of $G$.  Denote
\begin{equation} \label{def-a}
\ga: \text{ maximal abelian subspace of } 
\{\xi \in \gg \mid \theta(\xi) = -\xi\}.\\
\end{equation}
If $\xi \in \ga$ then $ad(\xi)$ is a semisimple linear transformation of
$\gg$ with all eigenvalues real.  Now, as usual, $\gg$ is the direct
sum of the joint eigenspaces (= restricted root spaces)
\begin{equation}
\gg^\gamma = \{\eta \in \gg \mid [\xi, \eta] = \gamma(\xi)\eta
\text{ for every } \xi \in \ga\} \text{ where } 
\gamma \text{ ranges over } \ga^*\ .
\end{equation}
The $\ga${\em --root system} of $\gg$ is $\Sigma(\gg,\ga) =
\{\gamma \in (\ga^* \setminus \{0\}) \mid \gg^\gamma \ne 0\}$.  Fix
\begin{equation}
\Sigma(\gg,\ga)^+: \text{ positive } \ga\text{\rm --root system of } \gg\ .
\end{equation}
Any two such systems are conjugate by the normalizer of $\ga$ in $K$.
\m

Then $\Sigma(\gg, \ga^+)$ specifies a nilpotent subalgebra
and a nilpotent subgroup, by
\begin{equation}
\gn = \sum_{\gamma \in \Sigma(\gg,\ga)^+}\ \gg^{-\gamma} \subset \gg
\ \ \text{ and } N \text{ is the analytic subgroup of } G \text{ for } \gn\ .
\end{equation}
The corresponding {\em minimal parabolic subalgebra} $\gp \subset \gg$
and {\em minimal parabolic subgroup} $P \subset G$ are given by
\begin{equation}
\gp \text{ is the normalizer of } \gn \text{ in } \gg \ \ 
\text{ and } \ \ P \text{ is the normalizer of } N \text{ in } G.
\end{equation}
Now denote
\begin{equation}
A: \text{ analytic subgroup of } G \text{ for } \ga \ \ \text{ and } \ \
M: \text{ centralizer of } A \text{ in } K.
\end{equation}
Then
\begin{equation}
\gp = \gm + \ga + \gn \ \ \text{ and } \ \ P = MAN \text{ with } 
MA = M \times A 
\end{equation}
and the corresponding {\em Iwasawa decompositions} are
\begin{equation}
\gg = \gk + \ga + \gn \ \ \text{ and } \ \ G = KAN.
\end{equation}

Both $M$ and $MA$ have the properties (\ref{group-class1}) and
(\ref{group-class2}).  Also, $M$ is compact modulo $Z_M(M^0)$, the
centralizer of $M^0$ in $M$.  We
write $\widehat{\phantom{X}}$ for unitary dual.  If $\xi \in \widehat{Z_{M^0}}$
we write $(\widehat{M^0})_\xi$ for the classes $[\eta^0] \in \widehat{M^0}$ 
such that $\eta^0|_{Z_{M^0}}$ is a multiple of $\xi$, and we write
$(\widehat{Z_M(M^0)})_\xi$ for the classes $[\chi] \in \widehat{Z_M(M^0)}$
such that $\chi|_{Z_{M^0}}$ is a multiple of $\xi$.  
\m

The extension of Cartan's highest weight theory appropriate for $M$ is

\begin{proposition} \label{m-structure}
{\em (Compare \cite[Proposition 1.1.3]{W1}.)}

{\rm 1.} $M = Z_M(M^0)M^0$.

{\rm 2.} Every irreducible representation of $M$ is
finite dimensional.

{\rm 3.} If $[\eta] \in \widehat{M}$
there exist unique $\xi \in \widehat{Z_{M^0}}$, 
$[\chi] \in (\widehat{Z_M(M^0)})_\xi$
and $[\eta^0] \in (\widehat{M^0})_\xi$ such that $[\eta] = [\chi \otimes \eta^0]$.

{\rm 4.} Let $\gt$ be a Cartan subalgebra of $\gm$, 
$\Sigma(\gm_{_\C}, \gt_{_\C})^+$ a
positive $\gt_{_\C}$--root system on $\gm_{_\C}$\,, and
$T^0 = \exp(\gt)$, so 
$\Lambda_\gm^+ = \{\nu \in \mathbf{i}\gt^* \mid e^\nu 
\text{ is well defined on } T^0 \text{ and } 
\langle \nu , \gamma \rangle \geqq 0 
\text{ for all } \gamma \in \Sigma(\gm_{_\C}, \gt_{_\C})^+\}$ is the set of 
dominant integral weights for $M^0$.  Then there is a bijection
$\nu \leftrightarrow [\eta_\nu]$ of 
$\Lambda_\gm^+$ onto $\widehat{M^0}$ given
by: $\nu$ is the highest weight of $\eta_\nu$.  Furthermore,
$[\eta_\nu] \in (\widehat{M^0})_\xi$ where
$\xi = e^\nu|_{Z_{M^0}}$.

{\rm 5.} $M = TM^0$ where $T$ is the Cartan subgroup 
$\{m \in M \mid Ad(m)\mu = \mu \text{, every } \mu \in \gt\}$ of $M$
that corresponds to the Cartan subalgebra $\gt$ of $\gm$.  
\end{proposition}

Define $\gh = \gt + \ga$.  It is a maximally split Cartan subalgebra of
$\gg$, and any two such Cartan subalgebras are $Ad(G^0)$--conjugate.
The positive root systems $\Sigma(\gg, \ga)^+$ and 
$\Sigma(\gm_{_\C}, \gt_{_\C})^+$ determine a
positive $\gh_{_\C}$--root system $\Sigma(\gg_{_\C},\gh_{_\C})^+$ for 
$\gg_{_\C}$ as follows.  A root $\gamma \in \Sigma(\gg_{_\C},\gh_{_\C})$
is positive if it
is nonzero and positive on $\ga$, or if it is zero on $\ga$ and positive
on $\gt_{_\C}$\,.   In other words,
\begin{equation}\label{coherent-root-order}
\begin{aligned}
&\Sigma(\gg, \ga)^+ = \{\gamma|_\ga \ \mid \ \gamma \in 
	\Sigma(\gg_{_\C},\gh_{_\C})^+ \text{ and }
		\gamma|_\ga \ne 0\} \ \text{ and } \\
&\Sigma(\gm_{_\C}, \gt_{_\C})^+ = \{\gamma|_\gt \ \mid \ \gamma \in 
	\Sigma(\gg_{_\C},\gh_{_\C})^+ \text{ and }
		\gamma|_\ga = 0\}.
\end{aligned}
\end{equation}
Now let
\begin{equation} \label{ps-data1}
[\eta_{\chi,\nu}] = [\chi \otimes \eta_\nu] \in \widehat{M} \ \ 
\text{and} \ \ \sigma \in \ga_{_\C}^*\ .
\end{equation}
That is equivalent to the datum
\begin{equation} \label{ps-data2}
\eta_{\chi,\nu,\sigma} \in \widehat{P} \text{ defined by }
\eta_{\chi,\nu,\sigma}(man) = e^\sigma(a)\eta_{\chi,\nu}(m) \text{ for }
m \in M, a \in A \text{ and } n \in N.
\end{equation}
Here $e^\sigma(\exp(\xi))$ means $e^{\sigma(\xi)}$ for $\xi \in \ga$.
In other words $e^\sigma(a)$ means $e^{\sigma(\log a)}$.
Also, we will write $V_{\chi,\nu, \sigma}$ for the representation space of
$\eta_{\chi,\nu,\sigma}$\, although as a vector space it is 
independent of $\sigma$.
\s

The corresponding {\em principal series} representation of $G$ is
\begin{equation} \label{ps_reps}
\pi_{\chi,\nu,\sigma} = \text{\rm Ind}_P^G(\eta_{\chi,\nu,\sigma}),
	\text{ induced representation.}
\end{equation}
Here one must be careful about the category in which one takes the induced 
representation.  For example, if $\cF$ is a smoothness class of functions
such as $C^k, 0 \leqq k \leqq \infty$, $C^\infty_c$ (test functions),
$C^{-\infty}$ (distributions),
$C^\omega$ (analytic), or $C^{-\omega}$ (hyperfunctions), one can take
$\pi_{\chi,\nu,\sigma}$ to be the natural representation (by translation of 
the variable) of $G$ on
\begin{equation} \label{smoothness-induced}
\cF(G,P:V_{\chi,\nu,\sigma}): \text{ all } f \in \cF(G : V_{\chi,\nu, \sigma}) 
\text{ with } f(gman) = e^{-\sigma}(a)\eta_{\chi,\nu}(m)^{-1}(f(g))
\end{equation}
for all $g \in G, m \in M, a \in A, \text{ and } n \in N$.
The representation is always given by the formula 
$\pi_{\chi,\nu,\sigma}(g)(f(g')) = f(g^{-1}g')$.
\m

One can also consider the analog of (\ref{smoothness-induced}) using
$K$--finite functions.  Those functions are $C^\omega$, and the 
representations spaces of the resulting $K$--finite
induced representations are the common underlying Harish--Chandra modules
for the representation spaces of the various smoothness classes 
(\ref{smoothness-induced}) of induced representations.
\m

Banach space representations, in particular unitary representations, are
more delicate.  We have to discuss this because we will have to keep track
of how they behave in a direct limit process.  
The modular function of $P$
is $\Delta_P(man) = e^{-2\rho_{\gg,\ga}}(a)$, where $2\rho_{\gg,\ga}(\xi)$
is the trace of $ad(\xi)|_\gn$ for $\xi \in \ga$, because we defined
$\gn$ to be the sum of the negative $\ga$--root spaces.
$G$ is unimodular, so $(\Delta_G/\Delta_P)(man) = \Delta_P^{-1}(man)
=  e^{2\rho_{\gg,\ga}}(a)$.
Let $\zeta$ be a norm--preserving representation of
$P$ on a Banach space $V_\zeta$ and let $1 \leqq p \leqq \infty$.
If $f \in C_c(G,P: \zeta \otimes \Delta_P^{-1/p})
= C_c(G,P: \zeta \otimes e^{(2/p) \rho_{\gg,\ga}})$ then 
$||f(\cdot)||_{V_\zeta}
\in C_c(G,P; \Delta_P^{-1/p})$, so the global norm
\begin{equation} \label{really-lp} \begin{array}{ccll}
||f||_p &=& \left (\int_{G/P} ||f(gP)||^p_{V_\zeta}
\, d\mu_{_{G/P}}(gP)\right )^{1/p} &\text{ for } p < \infty, \\
||f||_p &= &\text{\em ess sup\,}_{_{G/P}}||f(gP)||_{V_\zeta} &\text{ for }
p = \infty,
\end{array} \end{equation}
is well defined and invariant under the left translation action of
$G$.  Denote
\begin{equation*}
L_p(G,P:\zeta \otimes e^{(2/p) \rho_{\gg,\ga}}): \text{ Banach space completion of }
\left ( C_c(G,P: \zeta \otimes e^{(2/p) \rho_{\gg,\ga}}), ||\cdot||_p \right ).
\end{equation*}
Each $\pi_{\zeta \otimes e^{(2/p) \rho_{\gg,\ga}}}(g)$ extends by continuity 
from $ C_c(G,P: \zeta \otimes e^{(2/p) \rho_{\gg,\ga}})$ to a norm--preserving 
operator on $L_p(G,P:\zeta \otimes e^{(2/p) \rho_{\gg,\ga}})$ and that defines
a norm--preserving Banach representation of $G$ on 
$L_p(G,P,\zeta \otimes e^{(2/p) \rho_{\gg,\ga}})$.
If $\zeta$ is unitary then the global inner product
\begin{equation} \label{really-l2}
\langle f , f' \rangle = \int_{G/P} \langle f(gH), f'(gH)
        \rangle_{V_\zeta} \, d\mu_{_{G/P}}(gP) \text{ for } f, f' \in
        C_c(G,P: \zeta \otimes e^{\rho_{\gg,\ga}})
\end{equation}
is $G$--invariant, $L_2(G,P:\zeta \otimes  e^{\rho_{\gg,\ga}})$
is the Hilbert space completion of
$(C_c(G,P: \zeta \otimes  e^{\rho_{\gg,\ga}}), \langle \cdot , \cdot \rangle )$,
and $\pi_{\zeta \otimes  e^{\rho_{\gg,\ga}}}$ is a unitary representation of $G$.
Those unitary representations form the {\em unitary principal series} of $G$.
We translate the discussion to our terminology (\ref{ps_reps}) for principal
series representations as follows.

\begin{proposition}\label{unitary_criterion}
The principal series representation $\pi_{\chi,\nu,\sigma}$ extends by 
continuity from a representation of $G$ on $C_c(G,P:V_{\chi,\nu,\sigma})$ 
to a norm--preserving representation on $L_p(G,P:V_{\chi,\nu,\sigma})$ if
and only if $\sigma \in \mathbf{i}\ga^* + \frac{2}{p}\rho_{\gg,\ga}$\,.  In particular
it extends by continuity to a unitary representation of $G$ on
$L_2(G,P:V_{\chi,\nu,\sigma})$ if and only if $\sigma \in \mathbf{i}\ga^* + \rho_{\gg,\ga}$\,.
\end{proposition}

\begin{remark}{\em The restriction $\pi_{\chi,\nu,\sigma}|_K =
\text{\rm Ind}_M^K(\eta_{\chi,\nu})$, independent of $\sigma$.
Decompose $\sigma = \sigma' + \sigma''$ where $\sigma' \in \mathbf{i}\ga^* + \rho_{\gg,\ga}$ 
and $\sigma'' \in \ga$.  Then $\pi_{\chi,\nu,\sigma}|_K
= \pi_{\chi,\nu,\sigma'}|_K$, restriction of a unitary representation.
In other words, we may always view the underlying Harish--Chandra module
of $\pi_{\chi,\nu,\sigma}$ as a pre--Hilbert space.  This will be important
when we look at direct limit groups.}
\end{remark}

\section{Geometric Form of Principal Series Representations} \label{sec3}
\setcounter{equation}{0}

Let $G_{_\C}$ be a connected reductive complex Lie group for which $G$ is a
real form.  In other words there is a homomorphism $\varphi : G \to G_{_\C}$ 
with discrete kernel such that $d\varphi(\gg)$ is a real form of $\gg_{_\C}$\,.
If $Q$ is a parabolic subgroup of $G_{_\C}$\,, then we can view the complex
flag manifold $Z = G_{_\C}/Q$ as the set of all $G_{_\C}$--conjugates of
$Q$, say $Z \ni z \leftrightarrow Q_z \subset G_{_\C}$\,, because
$Q$ is its own normalizer in $G_{_\C}$\,.  Now we can view $Z$ as the 
set of all $Int(\gg_{_\C})$--conjugates of $\gq$ by
$Z \ni z \leftrightarrow \gq_z \subset \gg_{_\C}$\,.
\m

The condition (\ref{group-class1}) 
ensures that $G$ acts on $Z$ through $\varphi$ and conjugation.  In other
words $G$ acts on $Z$ through its adjoint action on $\gg_{_\C}$\,.  Thus
\begin{equation}\label{conjugation-action}
G \text{ acts on } Z \text{ by } g(z) = z' \text{ where }
\gq_{z'} = Ad(g)(\gq_z).
\end{equation}
This will be important when we construct direct limits of complex flag 
manifolds..
\m

Let $\Psi \subset \Sigma(\gm_{_\C}, \gt_{_\C})^+$ be any set of simple roots.
That defines a parabolic subalgebra $\gr = \gj_{_\C} + \gn_\gm$ in 
$\gm_{_\C}$\,,
with nilradical $\gn_\gm$ and Levi component $\gj_{_\C}$\,, where the reductive
algebra $\gj_{_\C}$ contains $\gt_{_\C}$ and has simple root system $\Psi$.
The corresponding parabolic subgroup of $M_{_\C}$ is $J_{_\C}N_\gm$, and
its $\varphi^{-1}$--image is a real form $J$ of $J_{_\C}$\,.
Note that $J = TJ^0$ where $T$ is the Cartan subgroup of $M$ corresponding
to $\gt$.
\m

Conversely to (\ref{coherent-root-order}) we extend roots of $\gm_{_\C}$
to roots of $\gm_{_\C} + \ga_{_\C}$ by zero on $\ga_{_\C}$ and obtain
$\Sigma(\gm_{_\C}, \gt_{_\C})^+ =
\Sigma(\gm_{_\C} + \ga_{_\C} , \gh_{_\C})^+ \subset
\Sigma(\gg_{_\C} , \gh_{_\C})^+$.  Every $\psi \in \Psi$ remains simple
for $\Sigma(\gg_{_\C} , \gh_{_\C})^+$.  Thus $\Psi$ also defines a
parabolic subalgebra $\gq = \gl_{_\C} + \gu$ in $\gg_{_\C}$\,,
with nilradical $\gu$ and Levi component $\gl_{_\C}$\,, where the reductive
algebra $\gl_{_\C}$ contains $\gh_{_\C}$ and has simple root system $\Psi$.
The corresponding parabolic subgroup of $G_{_\C}$ is $Q = L_{_\C}U$, and
its $\varphi^{-1}$--image is a real form $L$ of $L_{_\C}$\,.  Note that
$L = JA = AJ = ATJ^0 = HJ^0$ where $H = T \times A$ is the Cartan
subgroup of $G$ corresponding to $\gh$.
\m

As before, $Z$ is the complex flag manifold $G_{_\C}/Q$.  Let $z_0 = 1Q \in Z$.
Then the closed $G$--orbit in $Z$ is $F := G(z_0) = K(z_0)$.  We will realize
principal series representations on partially holomorphic vector bundles over 
$F$.  
\m

Define $S = M(z_0)$.  Note that $M$ acts on $Z$ as a compact group.  The 
basic properties of $S$, from \cite[Chapter 1]{W1}, are
(1) $S = M_{_\C}(z_0)$, so $S$ is a complex flag manifold, 
$S = M_{_\C}/R$ where $R = \{m \in M_{_\C} \mid m(z_0) = z_0\}$ is a
parabolic subgroup of $M_{_\C}/R$ with Lie algebra $\gr$ as described above,
(2) $S \cong M/J$ with $J$ as described above, and
(3) If $g, g' \in G$ and $gS$ meets $g'S$ then $gS = g'S$; and
$P = \{g \in G \mid gS = S\}$.  In fact, in the notation of \cite{W0} the
$gS$ are the holomorphic arc components of $F$.  Thus we have
\begin{proposition} \label{basic-fibration}
Define $\beta: F  \to G/P = \{gS \mid g \in G\}$ by $\beta(gz_0) = gS$. 
Then $\beta: F  \to G/P$ is a well defined $C^\omega$ fiber bundle with 
structure group $P$.  The fiber over $gP$ is $gS$, which is maximal
among complex submanifolds of $Z$ that are contained in $F$. 
\end{proposition}

Note that $T = Z_M(M^0)T^0$.  For $Z_M(M^0)$ centralizes $\gt$, thus is
contained in $T$, and if $t \in T$ then $Ad(t)|_{M^0}$ is an inner automorphism
of $M^0$ that fixes every $\mu \in \gt$, thus given by $Ad(t')|_{M^0}$
for some $t' \in T^0$, so $tT^0 \subset Z_M(M^0)T^0$.  As $T \subset J$
and $J \cap M^0 = J^0$ we have

\begin{lemma}\label{v-structure}{\em (Compare \cite[Proposition 1.1.3]{W1}.)} 

{\rm 1.} $J = Z_M(M^0)J^0$. 

{\rm 2.} If $[\zeta] \in \widehat{J}$
there exist unique $\xi \in \widehat{Z_{M^0}}$,
$[\chi] \in (\widehat{Z_M(M^0)})_\xi$ and $[\zeta^0] \in (\widehat{M^0})_\xi$ 
such that $[\zeta] = [\chi \otimes \zeta^0]$. 

{\rm 3.} Let $\Lambda_\gj^+ = \{\nu \in i\gt^* \mid e^\nu
\text{ is well defined on } T^0 \text{ and }
\langle \nu , \gamma \rangle \geqq 0
\text{ for all } \gamma \in \Sigma(\gj_{_\C}, \gt_{_\C})^+\}$, the set of
dominant integral weights for $J^0$.  Then there is a bijection
$\nu \leftrightarrow [\zeta_\nu^0]$ of $\Lambda_\gj^+$ onto $\widehat{J^0}$ 
given by: $\nu$ is the highest weight of $\zeta_\nu^0$.  Furthermore,
$[\zeta_\nu^0] \in (\widehat{J^0})_\xi$ where
$\xi = e^\nu|_{Z_{M^0}}$.
\end{lemma}

The set of $\gm$--nonsingular dominant integral weights for $J^0$ is
\begin{equation} \label{m-nonsingular}
(\Lambda_\gj^+)' = \{\nu \in \Lambda_\gj^+ \mid 
\langle \nu + \rho_{\gm,\gt} , \gamma \rangle \ne 0 \text{ for all } 
\gamma \in \Sigma(\gj_{_\C}, \gt_{_\C})\}
\end{equation}
where $\rho_{\gm,\gt}$ is half the sum of the roots in 
$\Sigma(\gm_{_\C}, \gt_{_\C})^+$.  If $\nu \in (\Lambda_\gj^+)'$ there is
a unique Weyl group element $w \in W(\gm,\gt)$ such that 
\begin{equation}
\widetilde{\nu} := w(\nu + \rho_{\gm,\gt}) - \rho_{\gm,\gt} \in \Lambda_\gm^+\ .
\end{equation}
We write $q(\nu)$ for the length $\ell(w)$ of that Weyl group element.
\m

Let $\nu \in (\Lambda_\gj^+)'$.  Let $\zeta_\nu$ denote the irreducible
representation of $J^0$ with highest weight $\nu$ as in Lemma
\ref{v-structure}.  Denote $\xi = e^\nu|_{Z_{M^0}}$ and choose
$[\chi] \in \widehat{Z_M(M^0)}_\xi$,.  Then
$[\zeta_{\chi,\nu}] = [\chi \otimes \zeta_\nu]$ is a well defined 
element of $\widehat{J}$.
Let $E_{\chi,\nu}$ denote the representation space.  Let 
$\sigma \in \ga_{_\C}^*$.  The isotropy subgroup of $G$ at $z_0$ is
$JAN$, and the representation $\zeta_{\chi,\nu,\sigma}(jan) =
e^\sigma(a)\zeta_{\chi,\nu}(b)$ of $JAN$ defines
\begin{equation} \label{def-bundle}
\E_{\chi,\nu,\sigma} \to F: \ G\text{--homogeneous vector bundle with
fiber } E_{\chi,\nu,\sigma} \text{ over } z_0 
\end{equation}
where $E_{\chi,\nu,\sigma}$ is the representation space $E_{\chi,\nu}$
of  $\zeta_{\chi,\nu,\sigma}$\,.
Note that $\E_{\chi,\nu,\sigma}|_{gS} \to gS$ is holomorphic, for every
fiber $gS$ of $F \to G/P$.  Initially one is tempted to define the 
corresponding sheaf as
$$
\begin{aligned}
&\cO_\gn(\E_{\chi,\nu,\sigma}) \to F: \text{ germs of $C^\infty$ functions }
	h : G \to E_{\chi,\nu,\sigma} \text{ such that } \\
&\phantom{XXXX} \text{\rm \phantom{i}(i) \phantom{X}} h(gjan) = 
	\zeta_{\chi,\nu,\sigma}(jan)^{-1}(h(g)) 
	\text{ for } g \in G \text{ and } jan \in JAN \\
&\phantom{XXXX} \text{\rm (ii) \phantom{X}} h(g;\xi) + 
	d\zeta_{\chi,\nu,\sigma}(\xi)h(g) = 0
	\text{ for } g \in G \text{ and } \xi \in (\gj + \ga + \gn)_{_\C}
\end{aligned}
$$
but that causes a number of technical problems, and it is better to use
hyperfunctions as in \cite{S3} and \cite{SW} to ensure that the differentials
in the cohomology of $\cO_\gn(\E_{\chi,\nu,\sigma}) \to F$ have closed
range.  The correct definition is
\begin{equation} \label{def-sheaf}
\begin{aligned}
&\cO_\gn(\E_{\chi,\nu,\sigma}) \to F: \text{ germs of $C^{-\omega}$ functions }
        h : G \to E_{\chi,\nu,\sigma} \text{ such that } \\
&\phantom{XXXX} \text{\rm \phantom{i}(i) \phantom{X}} h(gjan) =
        \zeta_{\chi,\nu,\sigma}(jan)^{-1}(h(g))
        \text{ for } g \in G \text{ and } jan \in JAN \\
&\phantom{XXXX} \text{\rm (ii) \phantom{X}} h(g;\xi) +
        d\zeta_{\chi,\nu,\sigma}(\xi)h(g) = 0
        \text{ for } g \in G \text{ and } \xi \in (\gj + \ga + \gn)_{_\C}\ .
\end{aligned}
\end{equation}
Apply the Bott--Borel--Weil Theorem to each 
$\E_{\chi,\nu,\sigma}|_{gS} \to gS$.  By elliptic regularity, use of 
hyperfunction coefficients results in the same cohomology as use of smooth
coefficients.  The result is

\begin{proposition} \label{ps-from-bbw}
{\rm (Compare \cite[Theorem 1.2.19]{W1}.)}
If $\nu \notin (\Lambda_\gj^+)'$ then $H^q(F;\cO_\gn(\E_{\chi,\nu,\sigma}))
= 0$ for every integer $q$.  If $\nu \in (\Lambda_\gj^+)'$,
then $H^q(F;\cO_\gn(\E_{\chi,\nu,\sigma}))= 0$ for $q \ne q(\nu)$, and
the natural action of $G$ on $H^{q(\nu)}(F;\cO_\gn(\E_{\chi,\nu,\sigma}))$ 
is infinitesimally equivalent $($same underlying Harish--Chandra module$)$ 
to the principal series representation $\pi_{\chi, \widetilde{\nu}, \sigma}$\,.
\end{proposition}

\section{Principal Series for Direct Limit Groups} \label{sec4}
\setcounter{equation}{0}

Consider a countable strict direct system 
$\{G_i,\phi_{k,i}\}_{i , k \in I}$
of reductive Lie groups.  Thus $I$ is 
a countable partially ordered set.  If $i, k \in I$ there
exists $\gamma \in I$ with $i \leqq \gamma$ and $k \leqq \gamma$.
Each $G_i$ is a reductive Lie
group.  If $i \leqq k$ then 
$\phi_{k,i}: G_i \to G_k$ is a continuous 
group homomorphism.  Then we have the direct limit group
$G = \varinjlim G_i$\,, with direct limit topology, and the 
$\phi_{k,i}$ specify continuous group homomorphisms 
$\phi_i : G_i \to G$.  The strictness condition is that the
homomorphisms $\phi_i$ are homeomorphisms onto their images.
So we may in fact view the $\phi_{k,i}$ as inclusions and view
$G$ as the union of the $G_i$\,,  and then the original topology
on each $G_i$ is the subspace topology.  In particular $G_i$
sits in $G$ as a closed (thus regularly embedded) submanifold.
\m

Countability of $I$ has two important consequences.  First, it guarantees
the existence of a $C^\omega$ (real analytic) Lie group structure on $G$.
See \cite{NRW1}, \cite{NRW2}, \cite{NRW3} and \cite{G}.  Second, 
it guarantees that $I$
either is finite or has a cofinal subset order--isomorphic to the positive
integers.  Whenever it is convenient we will replace $I$ by that subset;
this change in the defining direct system $\{G_i,\phi_{k,i}\}_{i , k \in I}$
has no effect on the direct limit group $G = \varinjlim G_i$\,.
\m

We always assume that every $G_i$ satisfies (\ref{group-class1})  
and (\ref{group-class2}).
\m

We have the corresponding strict direct system
$\{\gg_i, d\phi_{k,i}\}_{i , k \in I}$
of reductive Lie algebras, the direct limit algebra
$\gg = \varinjlim \gg_i$ with the direct limit topology,
and injective homomorphisms $d\phi_i : \gg_i \to \gg$ that are 
$C^\omega$ diffeomorphisms onto their images.  We also have the exponential
map $\exp: \gg \to G$, direct limit of the $\exp: \gg_i \to G_i$\,.
The $C^\omega$ Lie group structure on the limit group $G$ is specified
by the condition that $\exp: \gg \to G$ is a $C^\omega$ diffeomorphism from
a neighborhood on $0$ in $\gg$ onto a neighborhood of $1$ in $G$.
Again see \cite{NRW1}, \cite{NRW2}, \cite{NRW3} and \cite{G}.
\m

Consider a {\em compatible family of representations} 
$\{\pi_i , W_i,\psi_{k,i}\}_{i , k \in I}$
of $\{G_i,\phi_{k,i}\}$\,.  Thus $W_i$ is a
locally convex topological vector space (usually Hilbert or Fr\' echet),
$\{W_i,\psi_{k,i}\}_{i , k \in I}$ is a strict
direct system, $\pi_i$ is a continuous representation of $G_i$ on
$W_i$\,,  and
\begin{equation}
\text{ if } i \leqq k , g_i \in G_i \text{ and }
w_i \in W_i \text{ then }
\pi_k(\phi_{k,i}(g_i))(\psi_{k,i}(w_i))
= \psi_{k,i}(\pi_i(g_i)(w_i))
\end{equation}
That of course results in continuous injective linear maps 
$\psi_i : W_i \to W$ with closed image, where 
$W = \varinjlim W_k$\,.  We have the direct limit representation 
$\pi = \varinjlim \pi_i$ of $G$ on $W$ given by 
\begin{equation}
\text{ if } g = \phi_i(g_i) \in G \text{ and } 
w = \psi_i(w_i) \in W \text{ then }
\pi(g)w = \psi_i(\pi_i(g_i)(w_i)).
\end{equation}

We now examine the situation where the $\pi_i$ are principal series
representations of the $G_i$\,.  For that we need direct limits of
minimal parabolic subgroups.  
\m

As mentioned above we may assume $I = \{1, 2, 3, \ldots\}$ with the usual
order.  Then we recursively construct Cartan involutions $\theta_i$
of $\gg_i$ such that if $i \leqq k$ then 
$\theta_k|_{d\phi_{k,i}(\gg_i)}$ is $\theta_i$,
in other words $d\phi_{k,i}(\gk_i) = \gk_k \cap
d\phi_{k,i}(\gg_i)$.  We know that $\theta_i$ extends
uniquely to $G_i$ in such a way that its fixed point set $K_i$
has Lie algebra $\gk_i$\,, contains the kernel of the adjoint
representation of $G_i$\,, and meets every component of $G_i$\,.
Thus $K_i$ is the $G_i$--normalizer of $K_i^0 = K_i 
\cap G_i^0$\,.  Because of components, however, we must explicitly
assume that
\begin{equation} \label{k-lim}
\text{if } i \leqq k \text{ then } \phi_{k,i}(K_i)
\subset K_k \text{\,, so we have } K = \varinjlim K_i\ .
\end{equation}
While it is tempting to try to get around the assumption (\ref{k-lim})
by assuming that the $G_i$ are connected, we would still meet the
same problem with the groups $M_i$ indicated below.
\m

Now $d\phi_{k,i}$ maps the $(-1)$--eigenspace of $\theta_i$
into the $(-1)$--eigenspace of $\theta_k$, so we can recursively construct
maximal abelian subspaces $\ga_i \subset \{\xi \in \gg_i \mid
\theta_i(\xi) = -\xi\}$ as in (\ref{def-a}) such that
$d\phi_{k,i}(\ga_i) \subset \ga_k$ for $i \leqq k$.
Then the corresponding analytic subgroups satisfy
\begin{equation} \label{a-lim}
\text{if } i \leqq k \text{ then } \phi_{k,i}(A_i)
\subset A_k \text{\,, so we have } A = \varinjlim A_i\ .
\end{equation}

Note $d\phi_{k,i}(\ga_i) = \ga_k \cap
d\phi_{k,i}(\gg_i)$.  This allows us to recursively
construct a sequence of elements $\zeta_i \in \ga_i^*$ such that
$\langle \zeta_i , \alpha_i \rangle \ne 0$ for all 
$\alpha_i \in \Sigma(\gg_i, \ga_i)$ and $d\phi_{k,i}(\zeta_k) = \zeta_i$
for $i \leqq k$.  Taking roots where that inner product is positive
we have positive root systems $\Sigma(\gg_i, \ga_i)^+$
such that $d\phi_{k,i}$ maps every negative restricted root space
$\gg_i^{-\alpha_i}$ into $\gn_k :=
\sum_{\beta_k \in \Sigma(\gg_k,\ga_k)^+} \gg_k^{-\beta_k}$\,.
Again, the corresponding analytic subgroups satisfy
\begin{equation} \label{n-lim}
\text{if } i \leqq k \text{ then } \phi_{k,i}(N_i)
\subset N_k \text{\,, so we have } N = \varinjlim N_i\ .
\end{equation}

Essentially as before, let $M_i$ denote the centralizer 
$Z_{K_i}(A_i)$ of $A_i$ in $K_i$.  In general the
behavior of the $M_i$ (or even their identity components and Lie algebras)
under the $\phi_{k,i}$ is unclear.  Thus we explicitly
assume that
\begin{equation} \label{m-lim}
\text{if } i \leqq k \text{ then } \phi_{k,i}(M_i)
\subset M_k \text{\,, so we have } M = \varinjlim M_i\ .
\end{equation}
Now we put all this together.  Under the assumptions (\ref{k-lim}) and
(\ref{m-lim}) we have
\begin{equation}\label{coherent-iwasawa}
\begin{aligned}
\text{Iwasawa }&\text{decompositions } G_i = K_i A_i N_i
\text{ and minimal parabolics }
P_i = M_i A_i N_i  \\
&\text{ such that } \phi_{k,i} \text{ maps } K_i \to K_k\,,\,
M_i \to M_k\,,
A_i \to A_k \text{ and } N_i \to N_k \, .
\end{aligned}
\end{equation}
In particular,
\begin{equation}\label{lim-i}
\text{we have an Iwasawa decomposition } G = KAN
\end{equation}
and if $i \leqq k$ then $\phi_{k,i}(P_i) \subset P_k$ so
\begin{equation} \label{lim-p}
\text{we have the limit minimal parabolic } P = \varinjlim P_i = MAN\ .
\end{equation}
Here $G = KAN$ and $P = MAN$ mean
\begin{equation}
\begin{aligned}
&\text{\rm \;(i) } (k,a,n) \phantom{,} \mapsto \phantom{,} kan
\text{ is a } C^\omega \text{ diffeomorphism of } K \times A \times N
\text{ onto } G , \\
&\text{\rm (ii) } (m,a,n) \mapsto man
\text{ is a } C^\omega \text{ diffeomorphism of } M \times A \times N
\text{ onto } P .
\end{aligned}
\end{equation}

\begin{example}\label{ddl-examples} {\em The diagonal embedding direct 
limit groups described in \cite[Section 5]{NRW3}, and their extension to 
noncompact real forms, all satisfy (\ref{k-lim}) and (\ref{m-lim}), leading
to the limits and decompositions $G = KAN$ and $P = \varinjlim P_i
= MAN$ of (\ref{lim-i}) and (\ref{lim-p}).  For example, let $\{r_n\}$
and $\{s_n\}$ be sequences of integers $\geqq 0$ where $1 \leqq n < \infty$
and $r_n + s_n \geqq 1$.  Fix $k_1 > 1$, and recursively define
$k_{n+1} = r_nk_n+s_n$\,, define $G_n = SL(k_n;\R)$ and 
$\phi_{n+1,n}: G_n \to G_{n+1}$ by $\phi_{n+1,n}(g) = \text{ diag}
(g,\ldots ,g; 1, \ldots ,1)$ with $r_n$ of $g$'s and $s_n$ of $1$'s.
Here $K_n$ is the special orthogonal group $SO(k_n)$, $A_n$ consists of the
diagonal matrices of determinant $1$ with positive diagonal entries in 
$G_n$\,, $M_n$ consists of the diagonal matrices determinant $1$ with 
diagonal entries $\pm 1$ in $G_n$\,,
and $N_n$ consists of the lower triangular matrices in $G_n$ with
all diagonal entries equal to $1$.  The limit groups depend on the choice 
of sequences $\{r_n\}$ and $\{s_n\}$, and it is quite nontrivial to see
when pairs of sequences lead to isomorphic limits.
\hfill $\diamond$
}
\end{example}

In order to discuss representations of $M$ we need direct systems of 
Cartan subalgebras and appropriate root orders.  With $I = \{1, 2, \ldots \}$
we recursively construct
\begin{equation} \label{coherent-m-cartan}
\begin{aligned}
&\text{Cartan subalgebras } \gt_i \text{ in } \gm_i\ , 
d\phi_{k,i} \gt_i \subset \gt_k \text{, and positive 
systems }
  \Sigma(\gm_{i,\C}\ , \ \gt_{i, \C})^+ \\
&\text{such that }
d\phi_{k,i}\left ( {\sum}_{\alpha \in
\Sigma(\gm_{i,\C}\ ,\ \gt_{i, \C})^+} \gm_i^\alpha \right )
\subset {\sum}_{\beta \in \Sigma(\gm_{k,\C}\ ,\ \gt_{k,\C})^+}
\gm_k^\beta \text{ for } i \leqq k .
\end{aligned}
\end{equation}
Then $\gt = \varinjlim \gt_i$ is a Cartan subalgebra of $\gm$, the
root system $\Sigma(\gm_{_\C} , \gt_{_\C} ) = 
\varprojlim \Sigma(\gm_{i,\C}\ , \ \gt_{i, \C})$ (inverse
limit), and the positive system $\Sigma(\gm_{_\C} , \gt_{_\C} )^+ =
\varprojlim \Sigma(\gm_{i,\C}\ , \ \gt_{i, \C})^+$ is well
defined.  The Cartan subalgebra $\gt$ of $\gm$ defines a Cartan
subgroup $T = \{m \in M \mid Ad(m)\xi = \xi \text{ for all } \xi \in \gt\}$,
and $T^0 = T \cap M^0$ is the corresponding Cartan subgroup of $M^0$.
\m

Each $M_i/M_i^0$ is discrete, so $\varinjlim M_i^0$ is
connected, closed and open in $M$.  Thus $\varinjlim M_i^0 = M^0$
and $M/M^0$ is discrete.  The same considerations hold for $G$, $K$ and
$T$.  Also, since each $M_i = T_i M_i^0$ we have
$M = TM^0$\,.  Here note $T^0 = T \cap M^0$.  In the special case where each 
$\phi_{k,i}: Z_{M_i}(M_i^0) \to 
Z_{M_k}(M_k^0)$ we have $Z_M(M^0) 
= \varinjlim Z_{M_i}(M_i^0)$ and $M = Z_M(M^0)M^0$.
\m

A linear functional $\nu \in \gt_{_\C}^*$ is called {\em integral}
(or $\gm$--integral)
if $e^\nu$ is a well defined homomorphism $T^0 \to \C^\times$, in other
words if the pull--backs $\nu_i = 
\phi_i^*(\nu) \in \gt_{i, \C}^*$ are integral.  Here note $\nu =
\varprojlim \nu_i$\,.  The functional $\nu$ is called
{\em dominant} (or $\gm$--dominant) 
if $\langle \nu , a \rangle \geqq 0$ for every
$a \in \Sigma(\gm_{_\C} , \gt_{_\C} )^+$, in other words if $\nu_i$
is $\gm_i$--dominant for each $i$.   We use these notions for 
a small variation on the Mackey little--group method.

\begin{proposition} \label{rep-m}
Let $\nu \in \gt_{_\C}^*$ be a dominant integral linear functional.  It
determines an irreducible unitary representation $\eta_\nu$ of $M^0$
as follows.  Let $\eta_{i, \nu}$ denote
the irreducible unitary representation of $M_i^0$ with lowest weight
$- \nu_i = \phi_i^*(- \nu)$.  Choose a unit lowest weight
vector $v_{i,\nu}$ in the representation space $V_{i,\nu}$ of
$\eta_{i,\nu}$.  For $i \leqq k$ extend the Lie algebra monomorphism
$d\phi_{k,i} : \gg_i \hookrightarrow \gg_k$ as usual to an enveloping algebra 
monomorphism $\cU(\gg_i) \hookrightarrow \cU(\gg_k)$, which we also denote
$d\phi_{k,i}$\,, and define
$\psi_{k,i} : V_{i,\nu} \to V_{k;\nu}$ by 
$\psi_{k,i}(d\eta_{i,\nu}(\Xi_i)(v_{i,-\nu}))
= d\eta_{k,\nu}(d\phi_{k,i}(\Xi_i))(v_{k,-\nu})$
for $\Xi \in \cU(\gm_i)$. \hfill\newline
\phantom{X,} {\rm (1)} $\{\eta_{i,\nu}, V_{i,\nu}, \psi_{k,i}\}$ is a
compatible system of irreducible representations of
$\{M_i^0, \phi_{k,i}\}$, so
$\eta_\nu = \varinjlim \eta_{i,\nu}$ is a well defined irreducible
unitary representation of $M^0$, with representation space
$V_\nu = \varinjlim V_{i,\nu}$\,.  Further, $\eta_\nu$ is a lowest
weight representation with lowest weight $-\nu$, and
$v_{-\nu} := \varinjlim v_{i,-\nu}$ is a lowest weight unit vector. 

{\rm (2)} If $m \in M$ then $\eta_\nu\circ Ad(m)^{-1}$ is unitarily 
equivalent to $\eta_\nu$.  

{\rm (3)} Denote $\widehat{M_\nu}$: equivalence classes of irreducible 
unitary representations $\eta$ of $M$ such that $\eta|_{M^0}$ weakly
contains $\eta_\nu$ in the sense that the kernel of $d\eta$ on the
enveloping algebra $\cU(\gm)$ is contained in the kernel of $d\eta_\nu$\,.
Then $\widehat{M_\nu} = \{[\eta] \in \widehat{M} \mid [\eta|_{M^0}]
\text{ is a multiple of } \eta_\nu\}$.

{\rm (4)}
Let $[\eta_{\chi,\nu}] \in \widehat{M_\nu}$\,.  Let $V_{\chi,\nu}$ denote
its representation space.  Choose a subspace $V'_\nu \subset V_{\chi,\nu}$
on which $M^0$ acts by $\eta_\nu$\,, let $v \mapsto v'$ denote the
intertwining map of $V_\nu$ onto $V'_\nu$ and let $v'_{-\nu}$ be the image
of the lowest weight unit vector $v_{-\nu}$ of $\eta_\nu$\,.  Then the image
of $V_{i,\nu}$ in $V_{\chi,\nu}$ is 
$d\eta_{i,\nu}(\cU(\gm_i))(v'_{-\nu})$, and
$V'_\nu = \varinjlim d\eta_{i,\nu}(\cU(\gm_i))(v'_{-\nu})
= d\eta_\nu(\cU(\gm))(v'_{-\nu})$.
 
{\rm (5)} In the special case where $M = Z_M(M^0)M^0$, the set
$\widehat{M_\nu}$ consists of all $[\chi \otimes \eta_\nu]$
such that $\chi \in (\widehat{Z_M(M^0)})_\xi$ where
$\xi = e^{- \nu}|_{Z_{M^0}}$\,.
\s

\noindent
{\em {\sc Note.}  In general we write the elements of $\widehat{M_\nu}$ as 
$[\eta_{\chi,\nu}]$ where $\chi$ is just a parameter.  In the case of Statement
(5) the parameter $\chi$ is interpreted as an element of
$(\widehat{Z_M(M^0)})_\xi$\,.}
\end{proposition}

\noindent {\bf Proof.}  Statement (1) is satisfied by construction.
\s

For Statement (2)
let $m \in M$ and $\eta'_\nu = \eta_\nu\circ Ad(m)^{-1}$.  We view $M$ as 
the union of the $M_i$\,.  Then $m$ belongs to some $M_\delta$\,, 
hence to $M_i$ for $i \geqq \delta$.  Altering $m$ by an element
of $M_\delta^0$ we may assume that $m \in Z_{M_\delta}(M_\delta^0)$, so
$Ad(m)^*(\nu_\delta) = \nu_\delta$\,, and thus $\eta_{\delta,\nu}(v_\delta)$
is some multiple $c_\delta v_\delta$ of $v_\delta$\,.  Apply the 
enveloping algebra now to see that $v \mapsto c_\delta v$ intertwines
$\eta'_{\delta,\nu}$ with $\eta_{\delta,\nu}$\,.  The point here is that
we may replace $c_\delta$ by any other modular scalar, for example by $1$.
Now $v \mapsto v$ intertwines $\eta'_{i,\nu}$ with
$\eta_{i,\nu}$ for every $i \geqq \delta$,  and thus intertwines
$\eta'_\nu$ with $\eta_\nu$\,.
\s

For Statement (3) let $\cK$ denote the kernel of $d\eta$ on
$\cU(\gm)$ and let $\cK_\nu$ denote the kernel of $d\eta_\nu$\,.
If $\eta|_{M^0}$ is a multiple of $\eta_\nu$ then $\cK = \cK_\nu$\,.
Now let $\cK \subset \cK_\nu$\,.  Then the associative algebra 
$\cU(\gm)/\cK_\nu$ is a quotient of $\cU(\gm)/\cK$.  Remember that
$\eta_\nu$ is irreducible.
Since $M^0$ is connected and generated by $\exp(\gm)$ now $\eta_\nu$ 
is equivalent to a quotient representation of $\eta|_{M^0}$\,.  
By unitarity now $\eta_\nu$ is equivalent to a subrepresentation of 
$\eta|_{M^0}$\,.  Let $w$ be a cyclic unit vector for that irreducible
subrepresentation and let $W$ be a set of representatives of $M$ modulo $M^0$.
Then the representation
space of $\eta$ is generated by the $\eta(M)(\eta(x)w), x \in W$.
By Statement (2), the action of $M^0$ on the closed span of 
$\eta(M)(\eta(x)w)$ is
equivalent to $\eta_\nu$\,.  Thus $\eta|_{M^0}$ is a multiple of $\eta_\nu$\,.
\s

Statements (4) and (5) follow from (1) and (3).  \hfill $\square$
\m

Fix $[\eta_{\chi,\nu}] \in \widehat{M_\nu}$ as in Proposition \ref{rep-m}.
In the notation of Proposition \ref{rep-m}, identify $V_\nu$ with its
image $V'_\nu = d\eta_\nu(\cU(\gm))(v'_{-\nu})$ in $V_{\chi,\nu}$ and 
identify the lowest weight vector $v_{-\nu}$ of $\eta_\nu$ with its
image $v'_{-\nu}$ in $V_{\chi,\nu}$\,.  Let $V_{i,\chi,\nu}$
denote the closed span of $\eta_{\chi,\nu}(M_i)(v_{-\nu})$ and let 
$\eta_{i,\chi,\nu}$ denote the representation of $M_i$ on
$V_{i,\chi,\nu}$\,.  Unwinding the definitions one sees that
\begin{equation} \label{indprep-m}
\eta_{\chi,\nu} = \varinjlim \eta_{i,\chi,\nu}\ .
\end{equation}
Now let
\begin{equation} \label{ps-data1-lim}
\sigma \in \ga_{_\C}^* \text{ and } 
\sigma_i = \phi_i^*(\sigma) \in (\ga_i)_{_\C}^*\ .
\end{equation}
As in Section \ref{sec2} that is equivalent to the data
\begin{equation} \label{ps-data2-lim}
\begin{aligned}
&\eta_{\chi,\nu,\sigma} \in \widehat{P}: \  
\eta_{\chi,\nu,\sigma}(man) = e^\sigma(a)\eta_{\chi,\nu}(m) \text{ for }
m \in M, a \in A \text{ and } n \in N \text{, and } \\
&\eta_{i,\chi,\nu,\sigma} \in \widehat{P_i}: \ 
\eta_{i,\chi,\nu,\sigma}(man) = 
e^{\sigma_i}(a)\eta_{i,\chi,\nu}(m) \text{ for }
m \in M_i, a \in A_i \text{ and } n \in N_i
\end{aligned}
\end{equation}
Again as in Section \ref{sec2} we write $V_{\chi,\nu,\sigma}$ for the 
representation space 
of $\eta_{\chi,\nu,\sigma}$\,; as a vector space it just $V_{\chi,\nu}$\,.
Similarly we write $V_{i,\chi,\nu,\sigma}$ for the
representation space of $\eta_{i,\chi,\nu,\sigma}$\,.
\m
 
The {\em principal series} representation of $G$ defined by 
$[\eta_{\chi,\nu}] \in \widehat{M_\nu}$ and $\sigma \in \ga_{_\C}^*$ is
\begin{equation} \label{ps_reps-lim}
\pi_{\chi,\nu,\sigma} = \text{\rm Ind}_P^G(\eta_{\chi,\nu,\sigma}),
        \text{ induced representation.}
\end{equation}
This representation is always given by the formula
$\pi_{\chi,\nu,\sigma}(g)(f(g')) = f(g^{-1}g')$.  Of course we
also have the principal series representations 
$\pi_{i,\chi,\nu,\sigma} = \text{\rm Ind}_{P_i}^{G_i}
(\eta_{i,\chi,\nu,\sigma})$ of $G_i$\,.  
\m

The principal series representations $\pi_{\chi,\nu,\sigma}$ of
(\ref{ps_reps-lim}) has representation space that consists of an 
appropriate class of functions 
$f : G \to V_{\chi,\nu,\sigma}$ such that $f(gman) = e^{-\sigma(a)}
\eta_{\chi,\nu}(m)^{-1}\cdot f(g)$ for $g \in G$ and $man \in MAN = P$.
Here recall that $V_{\chi,\nu,\sigma}$ is the representation space
of $\eta_{\chi,\nu,\sigma}$.
View the representation space $V_{i,\chi,\nu,\sigma}$
of $\eta_{i,\chi,\nu,\sigma}$ as the closed $M_i$--invariant 
subspace of $V_{\chi,\nu,\sigma}$ generated by 
$\eta_{\chi,\nu,\sigma}(M_i)(v_{-\nu})$.  Then
the representation space of $\pi_{i,\chi,\nu,\sigma}$ is
the subspace of the representation space of $\pi_{\chi,\nu,\sigma}$\,,
given by $f(G_i) \subset V_{i,\chi,\nu,\sigma}$\,.  Since
$G$ is the union of the $G_i$ and $V_{\chi,\nu,\sigma}$ is the
union of the $V_{i,\chi,\nu,\sigma}$ we have proved
\begin{proposition} \label{dlim-ps_reps}
The principal series representations of a countable
strict direct limit are just the direct limits of principal series
representations.  Specifically,
$\pi_{\chi,\nu,\sigma} = \varinjlim \pi_{i,\chi,\nu,\sigma}$\ .
\end{proposition}

In dealing with principal series representations one must be very careful 
about the category in which he takes the induced representation.  
Smoothness categories such as $C^k, 0 \leqq k \leqq \infty$, 
$C^\infty_c$ (test functions), $C^{-\infty}$ (distributions),
$C^\omega$ (analytic), or $C^{-\omega}$ (hyperfunctions)
are still available for principal series representations of $G$, but 
anything involving integration over $G/P$ is
excluded.  We will get around this problem by constructing geometric 
realizations that provide $L_p$ versions of the principal series for $G$.

\section{Groups and Spaces for the Limit Principal Series} \label{sec5}
\setcounter{equation}{0}

The Iwasawa decompositions (\ref{coherent-iwasawa}), and the Cartan 
subalgebras $\gt_i \subset \gm_i$ and the positive
root systems of (\ref{coherent-m-cartan}) define 
\begin{equation} \label{coherent-g-cartan}
\begin{aligned}
&\text{Cartan subalgebras } \gh_i = \gt_i \times \ga_i 
\text{ in } \gg_i \text{ such that }
d\phi_{k,i} \text{ maps } \gh_i \to \gh_k \\
&\text{and positive root systems } \Sigma(\gg_{i,\C},\gh_{i, \C})^+
\text{ given by } (\ref{coherent-root-order}) \\
&\text{such that } 
d\phi_{k,i}\left ( {\sum}_{a \in 
\Sigma(\gg_{i,\C}\ ,\ \gh_{i, \C})^+)} \gg_{i,a}\right )
\subset {\sum}_{b \in \Sigma(\gg_{k,\C}\ ,\ \gh_{k,\C})^+}
\gg_{k,b} \text{ for } i \leqq k .
\end{aligned}
\end{equation}
Then $\gh = \varinjlim \gh_i$ is a Cartan subalgebra of $\gg$,
$\Sigma(\gg_{_\C}, \gh_{_\C}) := \varprojlim 
\Sigma(\gg_{i,\C}\ ,\ \gh_{i, \C})$ is its root system,
and $\Sigma(\gg_{_\C}, \gh_{_\C})^+  := \varprojlim
\Sigma(\gg_{i,\C}\ ,\ \gh_{i, \C})^+$ is a positive subsystem.
Further, we will need
\begin{equation}\label{coherent-cpx-parabolics}
\begin{aligned}
\gq_i = &\gl_{i,\C} + \gu_i \subset \gg_{i,\C}
\text{\,: \ parabolic subalgebras such that } \\
&\text{\rm \phantom{ii}(i) } \text{\rm the }\gq_i
\text{ are defined by sets } \Psi_i
\text{ of } \Sigma(\gm_{i,\C}\ ,\ \gt_{i, \C})^+
\text{--simple roots, as in \S\ref{sec3}}, \\
&\text{\rm \phantom{i}(ii) }
d\phi_{k,i}^{-1}(\gq_k) = \gq_i 
\end{aligned}
\end{equation}
Then $\gl_i := \gg_i \cap \gl_{i, \C}$ is a real 
form of $\gl_{i, \C}$\,, and 
$\gg_i \cap \gq_i = \gl_i +  \gn_i$\,,
as in Section \ref{sec3}.
\m

Let $G_{i,\C}$ denote the connected simply connected Lie group with
Lie algebra $\gg_{i,\C}$\,.  In general $G_i$ will not be a real 
form of $G_{i,\C}$ because in general $\gg_i \hookrightarrow
\gg_{i,\C}$ will not integrate to a homomorphism $G_i
\to G_{i,\C}$\,, but at least we have the connected complex 
simply connected group $G_\C = \varinjlim G_{i,\C}$ with Lie algebra
$\gg_{_\C} = \varinjlim \gg_{i,\C}$\,.
\m 

Let $Q_i$ be the parabolic subgroup with Lie algebra $\gq_i$\,,
and let $Z_i$ denote the complex flag manifold 
$G_{i,\C}/Q_i$\,.
Note that we would get the same $Z_i$ if we did this construction
starting with arbitrary complex Lie groups 
$\hskip 1 pt ' \hskip -1 pt G_{i,\C}$ for which the
$G_{i,\C}$ are the universal covering groups, in particular if we
started with any connected complex Lie group $\hskip 1 pt ' \hskip -1 pt
G_{i,\C}$ for
which $G_i$ is a real form.  For $Z_i$ can be identified as
the set of all $Int(\gg_{i.\C})$--conjugates of $\gq_i$
in  $\gg_{i,\C}$\,, with the action of $G_i$ given by 
conjugation as in (\ref{conjugation-action}).
\m

The reason for this indirection is that in general we cannot choose
a family of complex Lie groups $\hskip 1 pt ' \hskip -1 pt 
G_{i,\C}$\,, for which 
the $G_{i,\C}$ are the universal covering groups, such that
the $\hskip 1 pt ' \hskip -1 pt G_{i,\C}$ constitute a well defined 
direct system of complex
Lie groups and holomorphic homomorphisms $\hskip 1 pt ' \hskip -1 pt
\phi_{k,i}$ with
$d\phi_{k,i} = d\hskip 3 pt ' \hskip -1 pt \phi_{k,i}$\,.
\m

We now recall some structural information concerning the limit groups and limit
flags from \cite[Sections 1 \& 2]{NRW3}.
\s

The parabolic $Q_i = L_{i,\C} U_i$\,, semidirect product, where
$L_{i, \C}$ and $U_i$ are the respective complex analytic subgroups 
of $G_{i, \C}$ for $\gl_{i,\C}$ and $\gu_i$\,.
The direct systems $\{G_{i,\C}, \phi_{k,i}\}$ and 
$\{\gq_i, d\phi_{k,i}\}$ define direct systems
$\{L_{i,\C}, \phi_{k,i}\}$ and
$\{U_i, \phi_{k,i}\}$.
Let $Q = \varinjlim Q_i$\,, $L_\C = \varinjlim L_{i,\C}$\,, and
$U = \varinjlim U_i$\,.  Then $Q$, $L_\C$ and $U$ are closed complex 
analytic subgroups of $G$, and $Q = L_\C U$ semidirect product.
\m

We define a direct system $\{Z_i , \phi'_{k,i}\}$ by
$Z_i = \{Ad(g_i)\gq_i \mid g_i \in G_i\}$ and
$\phi'_{k,i}(z_i) = z_k$\,, where $z_i =
Ad(g_i)(\gq_i)$ gives $z_k = 
Ad(\phi_{k,i}(g_i))(\gq_k)$.  Then 
$\{Z_i , \phi'_{k,i}\}$ is a strict direct system of complex
manifolds and holomorphic maps, so the limit
$Z = \varinjlim Z_i$ is a complex manifold and the $\phi'_i :
Z_i \to Z$ are holomorphic injections with closed image.  The
$Z_i$ are complex homogeneous spaces $G_{i,\C}/Q_i$\,, and
it follows that the limit flag manifold $Z$ is a complex homogeneous
space $G_\C(z_0) = G_\C/Q$ where $z_0$ is the base point in $Z$, i.e.
$\phi'_i(z_{i,0}) = z_0$ for every $i$.  Further,
the action $G \times Z \to Z$ is holomorphic.
\m

Let $F_i$ denote the closed orbit $G_i(z_{i,0})
= K_i(z_{i,0})$ in $Z_i$\,.
Let $S_i = M_i(z_{i,0})$, complex flag manifold
$M_{i,\C}/R_i$ where $R_i$ is the parabolic subgroup of 
$M_{i,\C}$ for the set $\Psi_i$ of simple 
$(\gm_{i,\C},\gt_{i,\C})$--roots whose extension to
$\gh_{i,\C}$ defines $\gq_i$ as in (\ref{coherent-cpx-parabolics}).
We have $\gr_i = \gj_{i,\C} + \gn_{i,\gm}$\,, reductive part
and nilradical, and $\gj_{i,\C} = \gr_i \cap \gj_{i,\C}$ and
$\gn_\gm = \gr_i \cap \gn_{i,\gm}$\,.  Thus $R_i =
J_{i,\C}N_{i,\gm}$\,.  Up to finite covering, let $L_i$ denote
the real form $\phi_i^{-1}(L_{i,\C})$ of $L_{i,\C}$ and 
let $J_i$ denote the real form $\phi_i^{-1}(J_{i,\C})$ of 
$J_{i,\C}$\,.  Then $L_i = J_i \times A_i$\,.  
\m

Now go to the limit: $F = \varinjlim F_i$\,, 
$L = \varinjlim L_i$\,,
$S = \varinjlim S_i$\,,
$R = \varinjlim R_i$\,, and
$J_{_\C} = \varinjlim J_{i,\C}$\,.  Then 
$L = J \times A$, $G/LN = G/JAN \cong F \cong K/J$, and $S \cong M/J$.
Further, two translates $gS$ and $g'S$ either coincide or are disjoint, and 
$P = \{g \in G \mid gS = S\}$.  Thus we have a fibration exactly as in
Proposition \ref{basic-fibration}:

\begin{proposition} \label{basic-fibration-lim}
Define $k: F  \to G/P = \{gS \mid g \in G\}$ by $k(gz_0) = gS$.
Then $k: F  \to G/P$
is a well defined $C^\omega$ fiber bundle with structure group $P$, where
the fiber over $gS$ is the complex submanifold $gS$ of $Z$ that is contained
in $F$.
\end{proposition}
 
\section{Bundles and Sheaves for the Limit Principal Series} 
\label{sec6}
\setcounter{equation}{0}

Retain the notation of Section \ref{sec5}.
In order to construct a coherent family of homogeneous vector bundles
$\E_{i,\chi,\nu,\sigma} \to F_i$\,, we start with a coherent
family of representations, as in Proposition \ref{rep-m}. 
The proof of Proposition \ref{rep-b} just below, is the same as the
proof of Proposition \ref{rep-m}.

\begin{proposition} \label{rep-b}
Let $\nu \in \gt_{_\C}^*$ be a $\gj_{_\C}$--dominant integral linear 
functional.  It determines an irreducible unitary representation 
$\zeta_\nu$ of $J^0$ as follows.  Let $\zeta_{i, \nu}$ denote
the irreducible unitary representation of $J_i^0$ with lowest weight
$- \nu_i = \phi_i^*(- \nu)$.  Choose a unit lowest weight
vector $e_{i, -\nu}$ in the representation space $E_{i,\nu}$ of
$\zeta_{i,\nu}$.  For $i \leqq k$ define
$\psi_{k,i} : E_{i,\nu} \to E_{k,\nu}$ by \hfill\newline
$\psi_{k,i}(d\zeta_{i,\nu}(\Xi_i)(e_{i, -\nu}))
= d\zeta_{k,\nu}(d\phi_{k,i}(\Xi_i))(e_{k, -\nu})$
for $\Xi$ in the enveloping algebra $\cU(\gj_i)$. \hfill\newline
\phantom{X,} {\rm (1)} 
$\{\zeta_{i,\nu}, E_{i,\nu}, \psi_{k,i}\}$ is a
compatible system of irreducible representations of
$\{J_i^0, \phi_{k,i}\}$, so
$\zeta_\nu = \varinjlim \zeta_{i,\nu}$ is a well defined irreducible
unitary representation of $J^0$, with representation space
$E_\nu = \varinjlim V_{i,\nu}$\,.
 
{\rm (2)} If $j \in J$ then $\zeta_\nu\circ Ad(j)^{-1}$ is unitarily
equivalent to $\zeta_\nu$.
 
{\rm (3)} Denote $\widehat{J_\nu}$: equivalence classes of irreducible
unitary representations $\zeta$ of $J$ such that $\zeta|_{J^0}$ weakly
contains $\zeta_\nu$ in the sense that the kernel of $d\zeta$ on the
enveloping algebra $\cU(\gj)$ is contained in the kernel of $d\zeta_\nu$\,.
Then $\widehat{J_\nu} = \{[\zeta] \in \widehat{J} \mid [\zeta|_{J^0}]
\text{ is a multiple of } \zeta_\nu\}$.

{\rm (4)}
Let $[\zeta_{\chi,\nu}] \in \widehat{J_\nu}$\,.  Let $E_{\chi,\nu}$ denote
its representation space.  Choose a subspace $E'_\nu \subset E_{\chi,\nu}$
on which $J^0$ acts by $\zeta_\nu$\,, let $e \mapsto e'$ denote the
intertwining map of $E_\nu$ onto $E'_\nu$ and let $e'_{-\nu}$ be the image
of the lowest weight unit vector $e_{-\nu}$ of $\zeta_\nu$\,.  Then the image
of $E_{i,\nu}$ in $E_{\chi,\nu}$ is
$d\zeta_{i,\nu}(\cU(\gj_i))(e'_{-\nu})$, and
$E'_\nu = \varinjlim d\zeta_{i,\nu}(\cU(\gj_i))(v'_{-\nu})
= d\zeta_\nu(\cU(\gj))(e'_{-\nu})$.
 
{\rm (5)} In the special case where $J = Z_J(J^0)J^0$, the set
$\widehat{J_\nu}$ consists of all $[\chi \widehat{\otimes} \zeta_\nu]$
such that $\chi \in (\widehat{Z_J(J^0)})_\xi$ where
$\xi = e^{- \nu}|_{Z_{J^0}}$\,.
\s
 
\noindent
{\em {\sc Note.}  In general we write the elements of $\widehat{J_\nu}$ as
$[\zeta_{\chi,\nu}]$ where $\chi$ is just a parameter.  In the case of 
Statement (5) the parameter $\chi$ is interpreted as an element of
$(\widehat{Z_J(J^0)})_\xi$\,.}
\end{proposition}

Now let
\begin{equation} \label{vb-data1-lim}
[\zeta_{\chi,\nu}] \in \widehat{J_\nu} \text{ as in Proposition \ref{rep-b},
and let } \sigma \in \ga_{_\C}^*\ .
\end{equation}
As in Sections \ref{sec2} and \ref{sec4} that is equivalent to the datum
\begin{equation} \label{vb-data2-lim}
\zeta_{\chi,\nu,\sigma} \in \widehat{JAN} \text{ defined by }
\zeta_{\chi,\nu,\sigma}(jan) = e^\sigma(a)\zeta_{\chi,\nu}(b) \text{ for }
j \in J, a \in A \text{ and } n \in N.
\end{equation}
As in Sections \ref{sec2} and \ref{sec4} we write $E_{\chi,\nu, \sigma}$ 
for the representation space of $\zeta_{\chi,\nu,\sigma}$\,.  
Now we have
\begin{equation} \label{def-bundle-lim}
\E_{\chi,\nu,\sigma} \to F: \ G\text{--homogeneous vector bundle with
fiber } E_{\chi,\nu,\sigma} \text{ over } z_0
\end{equation}
as in (\ref{def-bundle}).  If $g \in G$ then $\E_{\chi,\nu,\sigma}|_{gS} \to
gS$ is a holomorphic vector bundle.
\m

Note that the limit $E_{\chi,\nu, \sigma} = 
\varinjlim E_{i, \chi,\nu, \sigma}$ 
where $E_{i, \chi,\nu, \sigma}$ is the subspace of $E_{\chi,\nu, \sigma}$
generated by $\zeta_{\chi,\nu,\sigma}(J_i A_i N_i)(e'_{-\nu})$.
Let $\E_{i, \chi,\nu, \sigma} \to F_i$ denote the associated
$G_i$--homogeneous vector bundle.  It is holomorphic over each 
$g_i S_i$\,.
The maps 
$$
\phi_{k,i} \times \psi_{k,i} : 
	G_i \times E_{i, \chi,\nu, \sigma} \to
	G_k \times E_{k, \chi,\nu, \sigma}
$$
induce bundle maps $(\phi_{k,i} \,, \psi_{k,i}):
        \E_{i, \chi,\nu, \sigma} \to \E_{k, \chi,\nu, \sigma}$\,.
These bundle maps form a coherent system and give us
\begin{equation} \label{bundle-dir-lim}
\E_{\chi,\nu,\sigma} = \varinjlim \E_{i, \chi,\nu, \sigma} \ .
\end{equation}  

Write $E_{\chi,\nu,\sigma}^*$ for the strong
topological dual of $E_{\chi,\nu,\sigma}$ and write
$\E_{\chi,\nu,\sigma}^* \to F$ for the associated homogeneous vector
bundle.  Again, the restricted bundles 
$\E_{\chi,\nu,\sigma}^*|_{gS} \to gS$ are holomorphic vector bundles, 
for every fiber $gS$ of $F \to G/P$, by Lemma 2.2 of \cite{NRW3}.
By elliptic regularity for hyperfunctions, Dolbeault cohomology is the
same for $C^\infty$ coefficients as for $C^{-\omega}$ coefficients.
Thus the corresponding sheaves are the
\begin{equation} \label{def-sheaf-lim}
\begin{aligned}
&\cO_\gn(\E_{\chi,\nu,\sigma}) \to F: \text{ germs of $C^{-\omega}$ functions }
        h : G \to E_{\chi,\nu,\sigma} \text{ such that } \\
&\phantom{XXXX} \text{\rm \phantom{i}(i) \phantom{X}} h(gjan) =
        \zeta_{\chi,\nu,\sigma}(jan)^{-1}(h(g))
        \text{ for } g \in G \text{ and } jan \in JAN \\
&\phantom{XXXX} \text{\rm (ii) \phantom{X}} h(g;\xi) +
        d\zeta_{\chi,\nu,\sigma}(\xi)h(g) = 0
        \text{ for } g \in G \text{ and } \xi \in (\gj + \ga + \gn)_{_\C}
\end{aligned}
\end{equation}
as in (\ref{def-sheaf}), and also the sheaf 
$\cO_\gn(\E_{\chi,\nu,\sigma}^*) \to F$ corresponding to the dual bundle.
These are the sheaves of germs of $C^{-\omega}$ sections that are holomorphic 
over the fibers $gS$ of $F \to G/P$.
\m

For simplicity of notation, we write $\cO_\gn(\E_{i,\chi,\nu,\sigma}) 
\to F_i$\,, instead of  $\cO_{\gn_i}(\E_{i,\chi,\nu,\sigma}) 
\to F_i$\,, for the sheaf over $F_i$ analogous to that of 
(\ref{def-sheaf-lim}).
\m

We recall the definition of the inverse limit sheaf
$\varprojlim \cO_\gn(\E_{i,\chi,\nu, \sigma}^*)$.  First, identify
$Z_i$ with $\phi_i(Z_i) \subset Z$, thus also identifying
$F_i$ with $\phi_i(F_i) \subset F$, and view
$\cO_\gn(\E_{i,\chi,\nu, \sigma}^*)$ as a sheaf over $F$ with stalk
$\{0\}$ over every point $z \notin F_i$\,.  The open subsets
of $F_i$ are the sets $U_i = U \cap F_i$ where $U$ is
open in $F$.  Let $\Gamma_i(U)$ denote the abelian group of sections of
$\cO_\gn(\E_{i,\chi,\nu, \sigma}^*)|_{U_i}$\,.  The $\Gamma_i(U)$
form a complete presheaf, corresponding to $\E_{i,\chi,\nu, \sigma}^*$\,.
Also, the abelian group $\Gamma(U)$ of sections of 
$\cO_\gn(\E_{\chi,\nu, \sigma}^*)|_U$ is the inverse limit, 
$\Gamma(U) = \varprojlim \Gamma_i(U)$\, corresponding to the inverse 
system given by restriction of sections and then extension by zero.
Also, the $\Gamma(U)$ form a
complete presheaf corresponding to $\cO_\gn(\E_{\chi,\nu, \sigma}^*)$.  Thus,
by definition,
\begin{equation} \label{proj-lim-sheaf}
\cO_\gn(\E_{\chi,\nu, \sigma}^*) = \varprojlim \cO_\gn(\E_{i,\chi,\nu, \sigma}^*).
\end{equation}

\begin{proposition}\label{coho-lim-is-lim-coho}
{\rm \sc (Compare \cite[Proposition 2.4]{NRW3})}
Let $q \geqq 0$.  Then there is a natural $G$--equivariant isomorphism
from the cohomology $H^q(F;\cO_\gn(\E_{\chi,\nu, \sigma}^*))$ of the inverse 
limit onto the inverse limit 
$\varprojlim H^q(F_i ; \cO_\gn(\E_{i, \chi,\nu, \sigma}^*))$
of the cohomologies.
\end{proposition}

\noindent {\bf Proof.}  Apply \cite[Chapter I, Theorem 4.5]{Ha} with the
global section functor $\Gamma$ in place of $T$ to see that our sheaf
cohomologies are the derived functors of $\Gamma$.  Our neighborhood bases
on $F$ and the $F_i$ are properly aligned, as described in the
above description of the  definition of the inverse limit sheaf, so that
we have a base $\cB$ for the topology of $F$ such that
each $\cB_i := \{U_i = U \cap F_i \mid U \in \cB\}$ forms
a base for the topology of $F_i$\,.   We can refine $\cB$ so that
the neighborhoods $U \in \cB$ have the following property.
If $U \in \cB$ and $g_i \in G_i$ such that $U \cap g_i S_i
\ne \emptyset$ then each
$U \cap g_i S_i$ is Stein, and $U_i$ is the product
of $(U \cap g_i S_i)$ with a cell.  Then, for every $U \in \cB$,

(a) the inverse system $\{\Gamma_i(U)\}$ is surjective, in other words
if $i \leqq k$ and $s_i \in \Gamma_i(U)$ then there
exists $s_k \in \Gamma_k(U)$ such that 
$s_i = s_k|_{U_i}$\,, and

(b) if $q > 0$ then 
$H^q(U_i, \cO_\gn(\E_{i, \chi,\nu, \sigma}^*)|_{U_i}) = 0$
for all $i$.
\m

The properties just noted are conditions (a) and (b) of 
\cite[Chapter I, Theorem 4.5]{Ha}.  Thus we have $G$--equivariant exact 
sequences
\begin{equation}
0 \to 
   {\varprojlim}^{(1)}H^{q-1}(F_i; \cO_\gn(\E_{i, \chi,\nu, \sigma}^*))
  \to H^q(F;\cO_\gn(\E_{\chi,\nu, \sigma}^*))
  \to \varprojlim H^q(F_i ; \cO_\gn(\E_{i, \chi,\nu, \sigma}^*)) \to 0
\end{equation}
where ${\varprojlim}^{(1)}$ denotes the first right derived functor of the
$\varprojlim$ functor.  The proof now is reduced to the proof that
${\varprojlim}^{(1)}H^{q-1}(F_i; \cO_\gn(\E_{i, \chi,\nu, \sigma}^*))
= 0$.  Following \cite[Chapter I, Theorem 4.3]{Ha} it suffices to check
the Mittag--Leffler condition
\begin{equation}\label{mittag--leffler}
\begin{aligned}
\text{for each } &i \text{ the filtration of } 
H^{q-1}(F_i; \cO_\gn(\E_{i, \chi,\nu, \sigma}^*)) \\
&\text{ by the }
H^{q-1}(F_k; \cO_\gn(\E_{k, \chi,\nu, \sigma}^*)) \text{ is
eventually constant.}
\end{aligned}
\end{equation}

Let $\eta_{i,\chi,\nu}$ denote the representation of $M_i$
on $V_i := 
H^{q-1}(S_i;\cO_\gn(\E_{i, \chi,\nu, \sigma}^*)|_{S_i})$.
If $\nu_i$ is $\gm_i$--singular then $V_i = 0$, so
$H^{q-1}(F_i; \cO_\gn(\E_{i, \chi,\nu, \sigma}^*)) = 0$ and the
Mittag--Leffler condition (\ref{mittag--leffler}) is trivially satisfied.
Now assume that $\nu_i$ is $\gm_i$--nonsingular.  From 
Proposition \ref{ps-from-bbw} (or see \cite[Theorem 1.2.19]{W1})
the action of $G_i$ on 
$H^{q-1}(F_i; \cO_\gn(\E_{i, \chi,\nu, \sigma}^*))$ is a certain
principal series representation.  Those representations have finite
composition series: see \cite[Theorem 4.4.4]{W1} for the unitary case,
and note that its proof suffices for the general case.  
The point there is that the infinitesimal
character and the $K_i$--restriction are fixed, and that forces
finiteness for the composition series.  Since each subspace in the
filtration $\{H^{q-1}(F_k; \cO_\gn(\E_{k, \chi,\nu, \sigma}^*)) \mid 
k \geqq i\}$ of 
$H^{q-1}(F_i; \cO_\gn(\E_{i, \chi,\nu, \sigma}^*))$ is an
$M_i$--submodule, there are only finitely many possible
composition factors, and (\ref{mittag--leffler}) is immediate.
That completes the proof of Proposition \ref{coho-lim-is-lim-coho}.
\hfill $\square$ 

\section{Geometric Realization of the Limit Principal Series} 
\label{sec7}
\setcounter{equation}{0}

In this section we establish the geometric realization of principal
series representations of direct limit groups, and look at some of
the consequences.  In effect we combine 
Propositions \ref{dlim-ps_reps}, \ref{basic-fibration-lim}, \ref{rep-b} and 
\ref{coho-lim-is-lim-coho}, and use ideas of Bott--Borel--Weil theory
from \cite{NRW3}.
\m

We first look at a limit construction for
principal series representations in the geometric style of the
limit Borel--Weil Theorem, where there is no problem of cohomology degree.

\begin{theorem}\label{limit-bw}
Let $\nu \in \gt_{_\C}^*$ be an $\gj_{_\C}$--dominant integral linear 
functional.  Let $\zeta_{\chi,\nu} \in \widehat{J_\nu}$ as in 
{\rm Proposition \ref{rep-b}.}  Let 
$\zeta_{\chi^*,\nu^*} = \zeta_{\chi,\nu}^*$\,,
the dual (contragredient) of $\zeta_{\chi,\nu}$\,.
For each $i$ suppose that $\nu_i := \phi_i^*(\nu)$ is
$\gm_{i, \C}$--dominant.  Let $\sigma \in \ga_{_\C}^*$\,.  Let $\sigma^*$
denote its complex conjugate, using conjugation of $\ga_{_\C}^*$ over
$\ga^*$.  Then the natural action of $G$ on 
$H^0(F;\cO_\gn(\E_{\chi,\nu, \sigma}^*))$ is infinitesimally equivalent 
to the principal series representation $\pi_{\chi^*,\nu^*,\sigma^*}
= \varprojlim \pi_{i, \chi^*,\nu^*,\sigma^*}$
of $G$, and its dual is infinitesimally equivalent
to the principal series representation 
$\pi_{\chi,\nu,\sigma} = \varinjlim \pi_{i,\chi,\nu,\sigma}$ of $G$.
\end{theorem}

\noindent {\bf Proof.}  Apply Proposition \ref{ps-from-bbw} to each
$\E_{i,\chi,\nu,\sigma} \to F_i$\,.  Since $\nu$ is dominant,
$\nu_i \in (\Lambda_{\gj_i}^+)'$, $q(\nu_i) = 0$,
and $\widetilde{\nu_i} = \nu_i$\,.  Thus the natural action
of $G_i$ on $H^0(F_i ; \cO_\gn(\E_{i,\chi,\nu,\sigma}))$
is the principal series representation $\pi_{i,\chi,\nu,\sigma}$\,.
\m

Note $\zeta_{\chi,\nu}^* = \zeta_{\chi^*,\nu^*}$
for some index $\chi^*$, and $\chi^*$ is in fact the dual of $\chi$ when
we are in the situation $J = Z_J(J^0)$ of Proposition \ref{rep-b}.  Also,
$e^{\sigma^*}$ is the dual of $e^\sigma$.  Thus the bundles 
$\E_{\chi,\nu,\sigma}$ and $\E_{\chi^*,\nu^*,\sigma^*}$ are dual, at least 
at the $K$--finite level.  Now $\pi_{\chi^*,\nu^*,\sigma^*}$
and $\pi_{\chi,\nu,\sigma}$ are dual, so the natural action of $G$ on
$H^0(F; \cO_\gn(\E^*_{\chi,\nu,\sigma}))$ is $\pi_{\chi^*,\nu^*,\sigma^*}$\,,
and the natural action of $G$ on
$H^0(F; \cO_\gn(\E_{\chi,\nu,\sigma}))$ is $\pi_{\chi,\nu,\sigma}$\,.
\m

Similarly $\pi_{i,\chi,\nu,\sigma}$ and 
$\pi_{i, \chi^*,\nu^*,\sigma^*}$ are dual, so the natural action of
$G_i$ on $H^0(F_i ; \cO_\gn(\E^*_{i,\chi,\nu,\sigma}))$
is $\pi_{i, \chi^*,\nu^*,\sigma^*}$\,.  Now  Proposition
\ref{coho-lim-is-lim-coho} says that $\pi_{\chi^*,\nu^*,\sigma^*}
= \varprojlim \pi_{i, \chi^*,\nu^*,\sigma^*}$\,, and thus also
$\pi_{\chi,\nu,\sigma} = \varinjlim \pi_{i,\chi,\nu,\sigma}$\,.
\hfill $\square$
\b

In order to extend Theorem \ref{limit-bw} to higher cohomology we face
the same problem as in \cite{NRW3}.  We have to find conditions under
which the cohomology degrees 
$$
q_i = q(\nu_i): = 
|\{ \gamma_i \in \Sigma(\gm_{i,\C},\gt_{i,\C})^+
\mid \langle \nu_i + \rho_{i, \gm,\gt} , \gamma_i \rangle < 0\}|
$$
remain constant as $i$ increases indefinitely.  So we recall some
definitions from \cite[Section 4]{NRW3}.
\m

Suppose that $\nu_i + \rho_{i,\gm,\gt}$ is nonsingular.  Then there
is a unique element $w_i$ in the Weyl group $W(\gm_i,\gt_i)$
that carries $\nu_i + \rho_{i,\gm,\gt}$ to a dominant weight, 
and $q_i = q_i(\nu_i)$ is the length $\ell(w_i)$.
\m

The {\em Weyl group} $W = W(\gm,\gt)$ is defined to be the group of all
$w|_\gt$ where $w$ is an automorphism of $\gm$ such that (i) $w(\gt) 
= \gt$ and (ii) for some index $i_0$ if $i \geqq i_0$ then
$w(d\phi_i(\gm_i)) = d\phi_i(\gm_i)$ and 
$w|_{d\phi_i(\gm_i)}$ is an inner automorphism of $\gm_i$\,.
\m

Our hypothesis (\ref{coherent-m-cartan}) amounts to a choice of Borel
subalgebra $\gb = \varinjlim \gb_i$ of $\gm$ such that $\gt_i \subset
\gb_i \subset \gr_i$ and $d\phi_{k,i}(\gb_i)
\subset \gb_k$\,, where $\gr_i = d\phi^{-1}_i(\gr)$.  
This choice determines the {\em finite Weyl group}
$W_F = W_F(\gm,\gb,\gt)$ consisting of all $w \in W$ such that 
$w(\gb) \cap \gb$ has finite codimension in $\gb$.  We define this
codimension to be the {\em length} $\ell(w)$.
\m

Let $w \in W(\gm,\gt)$.  Then we have the classically defined lengths
$\ell(w_i)$ relative to the positive root systems 
$\Sigma(\gm_{i,\C},
\gt_{i,\C})^+$.  If $w \in W_F(\gm,\gb,\gt)$ then there is an index
$i_0$\,, which in general depends on $w$, such that $\ell(w_i)
= \ell(w_k)$ for $k \geqq i \geqq i_0$\,, and this common
length is $\ell(w)$.  
\m

A linear functional $\nu \in \gt_{_\C}^*$ is {\em classically cohomologically 
finite} if there exist $w \in W_F(\gm,\gb,\gt)$ and $i_0$ as above,
and an integral linear functional $\widetilde{\nu} \in \gt_{_\C}^*$\,, with
the following property.  If $i \geqq i_0$ then 
$d\phi^*_i(\widetilde{\nu})$ is dominant relative to the positive
root system $\Sigma(\gm_{i,\C},\gt_{i,\C})^+$, and
$d\phi^*_i(\widetilde{\nu}) = w_i(\nu_i + \rho_{i,\gm,\gt})
- \rho_{i,\gm,\gt}$\,.  A linear functional $\nu \in \gt_{_\C}^*$ is {\em 
cohomologically finite} of {\em degree} $q_\nu$ if, whenever $i$ is 
sufficiently large, say $i \geqq i_0$\,, 
(i) $\nu_i + \rho_{i,\gm,\gt}$ is nonsingular and 
(ii) $q_i = q_\nu$ constant in $i$.  If $\nu$ is classically
cohomologically finite by means of $w \in W_F$ then it is cohomologically
finite of degree $\ell(w)$.  By contrast, there are cases where $\nu$ is
cohomologically finite of degree $q > 0$ while $W_F = \{1\}$, so $\nu$ is
not classically cohomologically finite.
\m

Drawing on \cite[Theorem 4.6]{NRW3} we now have a limit construction for
principal series representations in the geometric style of the 
Bott--Borel--Weil theorem, as follows.

\begin{theorem} \label{limit-bbw}
Let $\nu \in \gt_{_\C}^*$ be a $\gj_{_\C}$--dominant integral linear
functional.  Let $\sigma \in \ga_{_\C}^*$\,.
\s

\noindent {\rm 1.} If $\nu$ is not cohomologically finite then every
$H^q(F;\cO_\gn(\E^*_{\chi, \nu, \sigma})) = 0$.
\s

\noindent {\rm 2.} Assume that $\nu$ is cohomologically finite of degree
$q_\nu$\,.   Then
\s

{\rm (a)} $H^q(F;\cO_\gn(\E^*_{\chi, \nu, \sigma})) = 0$ for $q \ne q_\nu$\,,
and
\s

{\rm (b)} the natural action of $G$ on 
$H^{q_\nu}(F;\cO_\gn(\E^*_{\chi, \nu, \sigma}))$ is infinitesimally equivalent
to a principal series representation of the form 
$\pi_{\chi^*, \mu^*,\sigma^*}
= \varprojlim \pi_{i, \chi^*,\mu^*,\sigma^*}$\,,
and its dual is infinitesimally equivalent
to a principal series representation of the form
$\pi_{\chi,\mu,\sigma} = 
\varinjlim \pi_{i,\chi,\mu,\sigma}$\,.
\s

\noindent {\rm 3.} If further $\nu$ is classically cohomologically finite,
say by means of $w \in W_F$\,, then $q_\nu = \ell(w)$ and in {\rm (2)} we may
take $\mu = \widetilde{\nu}$, defined by $\mu_i = w_i(\nu_i
+ \rho_{i,\gm,\gt}) - \rho_{i,\gm,\gt}$ for $i$ sufficiently large.
\end{theorem}

\noindent {\bf Proof.}  Suppose that $\nu$ is not cohomologically finite.
Fix an integer $p \geqq 0$.  If $\nu_i + \rho_{i,\gm,\gt}$ is singular
then $H^p(S_i; \cO_\gn(\E^*_{\chi, \nu, \sigma}|_{S_i})) = 0$.  If
$\nu_i + \rho_{i,\gm,\gt}$ is nonsingular, then 
$H^p(S_i; \cO_\gn(\E^*_{\chi, \nu, \sigma}|_{S_i})) = 0$ unless
$q(\nu_i) = p$.  The $q(\nu_i)$ are increasing in $i$.
Since $\nu$ is not cohomologically finite, the $q(\nu_i)$ are unbounded.
Thus, $H^p(S_i; \cO_\gn(\E^*_{\chi, \nu, \sigma}|_{S_i}))$ becomes
$0$ and stays $0$ as $i$ increases.  Let $\eta^*_i$ denote the
representation of $M_i$ on 
$H^p(S_i; \cO_\gn(\E^*_{\chi, \nu, \sigma}|_{S_i}))$, and let
$\eta^*$ denote the representation of $M$ on 
$H^p(S; \cO_\gn(\E^*_{\chi, \nu, \sigma}|_S))$.  We have just seen that
$H^p(S; \cO_\gn(\E^*_{\chi, \nu, \sigma}|_S)) = 
\varprojlim H^p(S_i; \cO_\gn(\E^*_{\chi, \nu, \sigma}|_{S_i})) = 0$,
so the representation space of $\eta^*$ is $0$, and thus the
representation space $H^p(F;\cO_\gn(\E^*_{\chi, \nu, \sigma}))$ of
Ind$_{MAN}^G(\eta \otimes e^{\sigma^*})$ is zero.  That proves assertion (1).
\m

Assertion (2a) follows by an argument used for (1), and (2b) and (3) follow by
the argument of Theorem \ref{limit-bw}. \hfill $\square$
\m

\noindent Theorem \ref{limit-bbw} leaves us with two tasks:

1. find conditions on $\nu$ for cohomological finiteness, and

2. investigate boundedness and unitarity for the limit principal 
series representations.

\noindent The first is studied extensively in \cite{NRW3}, and we 
now turn to the second.

\section{Unitarity, ${\mathbf L}^p$ Boundedness, and Related Questions}
\label{sec8}
\setcounter{equation}{0}

According to Proposition \ref{unitary_criterion}, the $L_p$ condition 
for $\pi_{i,\chi,\nu,\sigma}$ is 
$\sigma_i \in \mathbf{i}\ga_i^* + \frac{2}{p}\rho_{i, \ga}$\,. 
So the $L_\infty$ condition is transparent: $\sigma_i \in \mathbf{i}\ga_i^*$ 
for all $i$ if and only if $\sigma \in \mathbf{i}\ga^*$. 
Now we set that case aside and suppose $1 \leqq p < \infty$.

\begin{lemma}\label{lp-cond}
Let $1 \leqq p < \infty$.  If $k \geqq i$ view 
$d\phi_{k,i} : \gg_i \to \gg_k$
as an inclusion $\gg_i \hookrightarrow \gg_k$\,.  Then the 
$\pi_{i,\chi,\nu,\sigma}$ satisfy
the $L_p$ condition for all $i \geqq i_0$ if and only if
{\rm (i)} 
$\sigma_{i_0} \in \mathbf{i}\ga_{i_0}^* + \frac{2}{p}\rho_{i_0, \ga}$
and {\rm (ii)} if $k \geqq i \geqq i_0$ then
$\rho_{k, \ga}|_{\ga_i} = \rho_{i,\ga}$\,.
In that case $\rho_\ga := \varprojlim \rho_{i,\ga} \in \ga^*$ is well defined.
\end{lemma}

\noindent {\bf Proof.}  The $\pi_{i,\chi,\nu,\sigma}$ satisfy
the $L_p$ condition for all $i \geqq i_0$ if and only if
$\sigma_i \in \mathbf{i}\ga_i^* + \frac{2}{p}\rho_{i,\ga}$ for
$i \geqq i_0$\,, in other words $Re\, \sigma_i = 
\frac{2}{p}\rho_{i,\ga}$\,.  
\s
If (i) and (ii) hold, it is obvious that the $\pi_{i,\chi,\nu,\sigma}$
satisfy the $L_p$ condition for all $i$ sufficiently large, say
$i \geqq i_0$\,.  Conversely,
suppose that $k \geqq i \geqq i_0$ and that 
$\pi_{\ell,\chi,\nu,\sigma}$ satisfies the $L_p$ condition for
$\ell = k, i, i_0$\,.  Then
$\frac{2}{p}\rho_{i,\ga} = \frac{2}{p}\rho_{k, \ga}|_{\ga_i}$\,,
in other words $\rho_{i,\ga} = \rho_{k,\ga}|_{\ga_i}$\,, in
addition to $\sigma_{i_0} \in \mathbf{i}\ga_{i_0}^* + 
\frac{2}{p}\rho_{i_0,\ga}$.  The last assertion follows.  \hfill $\square$
\m

Recall the structure theory for real parabolic subalgebras. 
Let $\Psi_i$ denote the set of simple roots of 
$\Sigma(\gg_i,\ga_i)^+$.  The $G_i^0$--conjugacy 
classes of (real) parabolic subalgebras of $\gg_i$ are in 
one to one correspondence $\Phi_i \leftrightarrow \gp_{i,\Phi}$ 
with the subsets $\Phi_i \subset \Psi_i$ by
\begin{equation}\label{real-parabolics}
\begin{aligned}
\gp_{i,\Phi} &= \gm_{i,\Phi} + \ga_{i,\Phi} + \gn_{i,\Phi}
		 \text{ where } \\
	& \ga_{i,\Phi} = \{\xi \in \ga_i \mid \psi_i(\xi) = 0 
		\text{ for all } \psi_i \in \Phi_i\}, \\
	& \gm_{i,\Phi} = \theta(\gm_{i,\Phi}) \text{ and } 
		\gm_{i,\Phi} \oplus \ga_{i,\Phi}
		\text{ is the centralizer of } \ga_{i,\Phi} \text{ in } 
		\gg_i\ , \\
	& \gn_{i,\Phi} \text{ is the sum of the negative 
		$\ga_i$--root spaces not in } \gm_{i,\Phi}\ .
\end{aligned}
\end{equation}
Here $\gn_{i,\Phi}$ is the nilradical, 
$\gm_{i,\Phi} \oplus \ga_{i,\Phi}$ is the Levi component
(reductive part), and $\gp_{i,\Phi}$ is the normalizer of 
$\gn_{i,\Phi}$ in $\gg_i$.  $\Phi_i$ is the simple root system
for $\gm_{i,\Phi} \oplus \ga_{i,\Phi}$\,.  The minimal parabolic
is the case $\Phi_i = \emptyset$.  The derived algebra 
$\gm'_{i,\Phi} = [\gm_{i,\Phi} , \gm_{i,\Phi}]$ is a \
maximal semisimple subalgebra of $\gp_{i,\Phi}$\,, and we refer to it
as the {\em semisimple component} of $\gp_{i,\Phi}$\,.
\m

If $\gamma \in \Sigma(\gg_i,\ga_i)$ we write
mult$(\gamma)$ for the multiplicity of $\gamma$ as an $\ga_i$--root,
in other words for the dimension $\dim \gg_i^\gamma$ of the root
space.  Thus $\rho_{i,\gg,\ga} = 
\sum_{\gamma \in \Sigma(\gg_i,\ga_i)^+}
\text{\rm mult}(\gamma)\gamma$. 
\m

The following lemma is standard in the context on non-restricted roots,
but we have not been able to find it in the literature, so we give a 
proof for the convenience of the reader.

\begin{lemma} \label{rho-prod}
If $\psi \in \Psi_i$ then 
$\frac{2\langle \rho_{i,\gg,\ga}, \psi \rangle}
{\langle \psi, \psi \rangle} = \text{\rm mult}(\psi) + 
2\,\text{\rm mult}(2\psi)$.
\end{lemma}

\noindent {\bf Proof.}  Let $w_\psi$ denote the Weyl group reflection
for the simple restricted root $\psi$.  Then 
$w_\psi\Sigma(\gg_i,\ga_i)^+ = \Sigma(\gg_i,\ga_i)^+
\setminus S(\psi)$ where $S(\psi)$ is $\{\psi\}$ if $2\psi$ is not a 
restricted root, $\{\psi, 2\psi\}$ if $2\psi$ is a restricted root.
Now compute
\begin{equation*}
\begin{aligned}
2 ( \rho_{i,\gg,\ga} -  &  \text{\rm mult}(\psi)\psi 
	- \text{\rm mult}(2\psi)2\psi ) 
		= 2 w_\psi(\rho_{i,\gg,\ga}) \\
& = {\sum}_{\gamma \in \Sigma(\gg_i,\ga_i)^+} w_\psi(\gamma) \\
& = {\sum}_{\gamma \in \Sigma(\gg_i,\ga_i)^+}
	\left ( \gamma - \frac{2\langle \gamma, \psi \rangle}
		{\langle \psi, \psi \rangle} \psi \right ) \\
& = 2\rho_{i,\gg,\ga} - \frac{2\langle 2\rho_{i,\gg,\ga},\psi\rangle}		{\langle \psi, \psi \rangle} \ .
\end{aligned}
\end{equation*}
Thus mult$(\psi) + 2$mult$(2\psi) = \frac{2\langle \rho_{i,\gg,\ga}, \psi \rangle}
{\langle \psi, \psi \rangle}$, as asserted.  
\hfill $\square$
\m

Now we are ready to look at Condition (ii) of Lemma \ref{lp-cond}.

\begin{proposition} \label{rho-restriction}
Let $\gg_i \subset \gg_k$\,, real semisimple Lie algebras.   Choose a
Cartan involution $\theta$ of $\gg_k$ that preserves $\gg_i$\,, let
$\ga_i$ be a maximal abelian subspace of 
$\{\xi \in \gg_i \mid \theta(\xi) = -\xi\}$, and enlarge $\ga_i$ 
to a maximal abelian subspace $\ga_k$ of
$\{\xi \in \gg_k \mid \theta(\xi) = -\xi\}$.  Suppose that
$\ga_k = \ga_i \oplus \ga_{k,i}$ where $\ga_{k,i}$ 
centralizes $\gg_i$\,, in other words that 
$\gg_i \oplus \ga_{k,i}$ is a subalgebra of $\gg_k$\,.
Then following conditions are equivalent.

{\rm 1.} The restriction $\rho_{k,\gg,\ga}|_{\ga_i} = 
\rho_{i,\gg,\ga}$\,.

{\rm 2.} $(\gg_i + \gm_k) \oplus \ga_{k,i}$ is 
the centralizer of $\ga_{k,i}$ in $\gg_k$\,.
 
{\rm 3.} Modulo $\gm_k$\,, the algebra $\gg_i$ is the semisimple 
component of a real parabolic subalgebra of $\gg_k$ that contains 
$\ga_k$\,. 
\end{proposition}

\noindent {\bf Proof.}  Assume (3).  Then there is a subset 
$\Phi \subset \Psi_k$ such that, modulo $\gm_k$\,, $\gg_i$
is the semisimple component $\gs$ of $\gp_{k,\Phi}$\,.  In particular
$\Phi$ is the simple root system for 
$\Sigma(\gg_i \oplus \ga_{k,i}, \ga_k)^+$,
so $\Sigma(\gg_i \oplus \ga_{k,i}, \ga_k) =
\Sigma(\gs \oplus \ga_{k,i}, \ga_k)$,
and the multiplicities
mult$_{\gg_i}(\gamma)$ = mult$_\gs(\gamma)$ for every root
$\gamma \in \Sigma(\gs \oplus \ga_{k,i}, \ga_k)$.
Thus $\rho_{i,\gg,\ga} = \rho_{\gs, \ga_i}$\,.
But Lemma \ref{rho-prod} shows that 
$\langle \rho_{k,\gg,\ga}, \varphi \rangle =
\langle \rho_{\gs \oplus \ga_{k,i},\ga},\varphi \rangle$ for
every $\varphi \in \Phi$, so $\rho_{k,\gg,\ga}|_{\ga_i} =
\rho_{\gs, \ga_i}$\,.  That proves (1).
\m

Assume (1).  Denote $\gr = (\gg_i + \gm_k)
\oplus \ga_{k,i}$\,.  We have not yet proved that $\gr$ is an
algebra, but we do have $\rho_{\gr,\ga_k} := \frac{1}{2}
\sum_{\gamma \in \Sigma(\gr,\ga_k)} \dim (\gr\cap\gg_k^\gamma)\gamma$,
and 
$\rho_{\gr,\ga_k} = \rho_{\gg_i \oplus \ga_{k,i},\ga_k}$
by definition of $\gr$.
\s

Let $\gz$ denote the centralizer of
$\ga_{k,i}$ in $\gg_k$\,.  Then
$\rho_{k,\gg,\ga}|_{\ga_i} = \rho_{\gz,\ga_k}|_{\ga_i}$\,.
Using assumption (1) now $\rho_{\gr,\ga_k} = \rho_{\gz,\ga_k}$\,.
By construction of $\gr$ and of $\gz$, if $\gamma \in \Sigma(\gr,\ga_k)^+$
then $\gamma \in \Sigma(\gz,\ga_k)^+$ and its multiplicities satisfy
mult$_\gr(\gamma) \leqq \text{ mult}_\gz(\gamma)$.  As 
$\rho_{\gr,\ga_k} = \rho_{\gz,\ga_k}$ now
$\sum_{\gamma \in \Sigma(\gz,\ga_k)^+} [\text{ mult}_\gz(\gamma)
- \text{ mult}_\gr(\gamma)]\gamma = 0$.  Take inner product with
$\rho_{\gz,\ga_k}|_{\ga_i}$.  Since 
 each $\langle \rho_{\gz,\ga_k}|_{\ga_i} , \gamma \rangle > 0$ and
 each $\text{mult}_\gz(\gamma) \geqq \text{ mult}_\gr(\gamma)$
it follows that
$\text{mult}_\gz(\gamma) = \text{ mult}_\gr(\gamma)$.
That proves $\gr = \gz$, which is the assertion of (2).
\m

Assume (2).  Then $\gz = (\gg_i + \gm_k) \oplus \ga_{k,i}$
is the reductive component of a parabolic subalgebra of $\gg_k$ and
the corresponding semisimple component is $[\gz,\gz] = 
[\gg_i + \gm_k , \gg_i + \gm_k]$.  If
$\gamma \in \Sigma(\gg_i \oplus \ga_{k,i}, \ga_k)$ then
$(\gg_i \oplus \ga_{k,i})^\gamma = \gg_k^\gamma$, so
$[\gm_k, \gg_i] \subset \gg_i$\,.  Now
$[\gz,\gz] = \gg_i + [\gm_k , \gm_k]$.  That proves (3),
completing the proof of the Proposition.
\hfill $\square$

\begin{corollary}\label{red-rho-parab}
Let $\gg_i$ be the semisimple component of a real parabolic
subalgebra of $\gg_k$ that contains $\ga_k$\,.  Then the
restriction $\rho_{k,\gg,\ga}|_{\ga_i} = \rho_{i,\gg,\ga}$\,.
\end{corollary}

\begin{definition} \label{def-w-parab}  {\em The strict direct system 
$\{G_i, \phi_{k,i}\}$ of reductive Lie groups 
is {\em weakly parabolic} if for
every pair $k \geqq i$ the subalgebra
$d\phi_{k,i}(\gg'_i) \hookrightarrow \gg'_k$ satisfies
the conditions of Proposition \ref{rho-restriction}, where $\gg'_\gamma$ 
denotes the derived algebra $[\gg'_\gamma, \gg'_\gamma]$  It is
{\em parabolic} if for every pair $k \geqq i$ the subalgebra
$d\phi_{k,i}(\gg'_i)$ is the semisimple component of a
real parabolic subalgebra of $\gg_k$\,.}
\end{definition}

\begin{remark} {\em The condition that $\{G_i, \phi_{k,i}\}$
be weakly parabolic is slightly less restrictive than the corresponding
condition (7.1) of ``coherent root orderings'' in \cite{NRW3}.  The context 
and applications are different, but the core idea is similar.}
\end{remark}

Now we come to the main result of this section:

\begin{theorem}\label{lp-limit}
Suppose that the direct system $\{G_i, \phi_{k,i}\}$ is weakly
parabolic.  Let $\nu \in \gt_{_\C}^*$ be a 
$\gj_{_\C}$--dominant integral linear functional
that is cohomologically finite of degree $q_\nu$\,.   
Let $\sigma \in \ga_{_\C}^*$\,.  Recall the principal series
representation $\pi_{\chi^*, \nu^*,\sigma^*}
= \varprojlim \pi_{i, \chi^*,\nu^*,\sigma^*}$ of $G$
on $H^{q_\nu}(F;\cO_\gn(\E^*_{\chi, \nu, \sigma}))$ and its dual
$\pi_{\chi,\nu,\sigma} = \varinjlim \pi_{i,\chi,\nu,\sigma}$\,.  
Suppose that $\sigma \in \mathbf{i}\ga^* + \frac{2}{p}\rho$, or
equivalently that some
$\sigma_i \in \mathbf{i}\ga_i^* + \frac{2}{p}\rho_{i, \ga}$\,,
where $1 \leqq p \leqq \infty$.  Then $\pi_{\chi, \mu,\sigma}$
is infinitesimally equivalent to a Banach space representation on
$\varinjlim L_p(G_i,P_i: V_{i,\chi, \mu,\sigma})$.  In 
particular, if $\sigma_{i_0} \in \mathbf{i}\ga_{i_0}^* + \rho_{i_0, \ga}$
then $\pi_{\chi, \mu,\sigma}$ is infinitesimally equivalent to a unitary
representation on 
$\varinjlim L_2(G_i,P_i: V_{i,\chi, \mu,\sigma})$.
\end{theorem}

\noindent {\bf Proof.} Combine Theorem \ref{limit-bbw} with 
Lemma \ref{lp-cond}, Proposition \ref{rho-restriction} and
Definition \ref{def-w-parab}.  \hfill $\square$
\m

In Theorem \ref{lp-limit} it would be better to derive the $L_p$ norm 
directly from the limit bundle $\E_{\chi, \mu,\sigma} \to F$.  We do this by
using a partially holomorphic cohomology space, as in \cite{W1}. 
The fibers $gS_i$ of $F_i \to G_i/P_i$ are compact,
so any cohomology class $c_{gS_i}  \in 
H^q(gS_i, \cO_\gn(\E_{\chi, \mu,\sigma}|_{gS_i}))$ is represented
by a harmonic $\E_{\chi, \mu,\sigma}|_{gS_i}$--valued $(0,q)$--form
$\omega_{gS_i}$\,, and $\omega_{gS_i}$ has a well--defined
$L_p$ norm $||\omega_{gS_i}||_p = 
( \int_{M_i} ||\omega_{gS_i}(m)||^p dm )^{1/p}$.
We now look at the Banach spaces
\begin{equation} \label{bpart-holo}
\begin{aligned}
\cB^q_{i,p}(F_i;\E_{i,\chi, \mu,\sigma}):\ &
	\text{measurable } \E_{i,\chi, \mu,\sigma}
	\text{--valued } (0,q) \text{--forms } \omega \text{ on } F_i
	\text{ such that } \\
& \omega|_{gS_i} \text{ is harmonic in the sense of Hodge and Kodaira,}\\
& ||\omega|_{gS_i}||_p \text{ is a measurable function of }
	gS_i \in G_i/P_i = K_i/M_i \ , \\
& \text{and } 
	\int_{K_i/M_i}  (||\omega_{gS_i}||_p)^p dk < \infty
\end{aligned}
\end{equation}
where 
\begin{equation} \label{banach-norm}
||\omega||_p = 
\left ( \int_{K_i/M_i}  (||\omega_{gS_i}||_p)^p dk\right )^{1/p}
= \left (\int_{K_i/J_i} (||\omega(kJ_i)||)^p dk\right )^{1/p}.
\end{equation}  
For $p = 2$ the norm is given by the inner product
\begin{equation}\label{inner-product}
\langle \omega , \omega' \rangle = {\int}_{K_i/M_i}\left (
{\int}_{M_i/J_i}\omega(kmz_0)\bar{\wedge}\#\omega(kmz_0)dm\right )dk.
\end{equation}
There $\#\omega$ is the $\E^*_{i,\chi, \mu,\sigma}$--valued 
$(s,s-q)$--form,
$s = \dim S_i$, which along $kS_i$ is the Hodge--Kodaira 
orthogonal of $\omega$, and $\bar{\wedge}$ is exterior product followed by
pairing of $\E_{i,\chi, \mu,\sigma}$ with $\E^*_{i,\chi, \mu,\sigma}$\,.
That gives us a Hilbert space
\begin{equation} \label{hpart-holo}
\cH^q_{i,2}(F_i;\E_{i,\chi, \mu,\sigma}) =
(\cB^q_{i,2}(F_i;\E_{i,\chi, \mu,\sigma}), 
\langle \cdot , \cdot \rangle ).
\end{equation}

Now, if we stay with a cofinal weakly parabolic subsystem of
$\{G_i,\phi_{k,i}\}$ as in Theorem \ref{lp-limit}, we have
Banach space representations
\begin{equation}\label{b-reps}
\begin{aligned}
&\pi_{\chi^*, \mu^*,\sigma^*} \text{ on }
  \cB^q_p(F;\E_{\chi^*, \mu^*,\sigma^*}) = 
  \varprojlim \cB^q_{i,p}(F_i;\E_{i,\chi^*, \mu^*,\sigma^*})\\
&\text{and } \pi_{\chi,\nu,\sigma} \text{ on }
\cB^q_{p'}(F;\E_{\chi, \mu, \sigma}) =
\varinjlim \cB^q_{i,p'}(F_i; \E_{i,\chi,\mu,\sigma})
\end{aligned}
\end{equation}
where $\frac{1}{p} + \frac{1}{p'} = 1$ for $1 < p < \infty$.  
In the case $p = 2$ we have unitary representations
\begin{equation}\label{u-reps}
\begin{aligned}
&\pi_{\chi^*, \mu^*,\sigma^*} \text{ on }
  \cH^q_2(F;\E_{\chi^*, \mu^*,\sigma^*}) =
  \varprojlim \cH^q_{i,2}(F_i;\E_{i,\chi^*, \mu^*,\sigma^*})\\
&\text{and } \pi_{\chi,\nu,\sigma} \text{ on }
\cH^q_2(F;\E_{\chi, \mu, \sigma}) =
\varinjlim \cH^q_{i,2}(F_i; \E_{i,\chi,\mu,\sigma}).
\end{aligned}
\end{equation}
Here note that $\varinjlim$ and $\varprojlim$ are the same in the Hilbert
space category.

\section{Diagonal Embedding Direct Limits}
\label{sec9}
\setcounter{equation}{0}

In this section we study an important class of direct limit groups that
includes those obtained from weakly parabolic direct systems.  These
diagonal embedding direct limits were introduced on the complex Lie algebra
level (see, for example, \cite{BS1}, \cite{BS2}, \cite{Ba}, \cite{BaZ}, 
\cite{YD} and \cite{Z}).  This topic is now plays a central role in the
theory of locally finite Lie algebras.  The idea was 
somewhat extended and applied on both the compact and the complex group 
level in \cite[Section 5]{NRW3}, and that is our starting point.
\m

{\sc Linear groups.}  
We consider limits of real, complex and quaternionic special linear groups.
Fix sequences ${\mathbf r} = \{r_n\}_{n \geqq 1}$,
${\mathbf s} = \{s_n\}_{n \geqq 1}$ and
${\mathbf t} = \{t_n\}_{n \geqq 1}$ of non--negative integers with
all $r_n + s_n > 0$.  Start with  $d_0 > 0$ and recursively 
define $d_{n+1} = d_n(r_{n+1} + s_{n+1}) + t_{n+1}$\,.  Let $\F$ be one of
$\R$ (real), $\C$ (complex) or $\H$ (quaternions) and define
$G_n = SL(d_n;\F)$.  Let $\delta$ denote the outer automorphism of
$G_n$ given by
\begin{equation} \label{outer}
\delta(g) = J(g^t)^{-1}J^{-1} \text{ where } J =
	\left ( \begin{smallmatrix} 0 & 0 & \hdots & 0 & 1 \\ 
		0 & 0 & \hdots & 1 & 0 \\ 
		\hdots & \hdots & \hdots & \hdots & \hdots & \\
        	0 & 1 & \hdots & 0 & 0 \\ 1 & 0 & \hdots & 0 & 0 
	\end{smallmatrix} \right ) .
\end{equation}
(The point of $J$ here is that $\delta$, as defined, preserves the
standard positive root system.) Then we have strict
direct systems $\{G_m,\phi_{n,m}\}_{n \geqq m \geqq 0}$ given by
\begin{equation}\label{ddl-lineargroups}
\phi_{n+1,n}: G_n \to G_{n+1} \text{ by }
\phi_{n+1,n}(g) \ = \diag\{g, \dots , g; \delta(g) , \dots , \delta(g) ;
1, \dots , 1\}
\end{equation}
with $r_{n+1}$ blocks $g$, with $s_{n+1}$ blocks $\delta(g)$, and with 
$t_{n+1}$ entries $1$.  That gives us 
\begin{equation} \label{ddl-lineargroups-lim}
G = SL_{{\mathbf r},{\mathbf s},{\mathbf t}}(\infty; \F) := 
\varinjlim \{G_m,\phi_{n,m}\} \text{ for the given } {\mathbf r}, {\mathbf s}
\text{ and } {\mathbf t}.
\end{equation}
Thus we have $SL_{{\mathbf r},{\mathbf s},{\mathbf t}}(\infty; \R)$,
$SL_{{\mathbf r},{\mathbf s},{\mathbf t}}(\infty; \C)$ and
$SL_{{\mathbf r},{\mathbf s},{\mathbf t}}(\infty; \H)$.
Of course the situation is exactly the same to construct infinite general
linear groups $GL_{{\mathbf r},{\mathbf s},{\mathbf t}}(\infty; \R)$
and $GL_{{\mathbf r},{\mathbf s},{\mathbf t}}(\infty; \C)$.
\m

{\sc Unitary groups.} 
We consider limits of real, complex and quaternionic unitary groups.  
Here $SU(p,q;\R)$ denotes the special orthogonal group $SO(p,q)$ for a
nondegenerate bilinear form of signature $(p,q)$,
$SU(p,q;\C)$ denotes the usual complex special unitary $SU(p,q)$ for a
nondegenerate hermitian form of signature $(p,q)$, and $SU(p,q;\H)$
is the quaternionic special unitary group for a nondegenerate hermitian
form signature $(p,q)$.  In each case we write the form as
$b(z,w) = \sum_{1\leqq 1\leqq p} \overline{w_i}z_i 
- \sum_{1\leqq 1\leqq q} \overline{w_{p+i}}z_{p+i}$\,, reflecting the fact
that we view $\F^{p+q}$ as a right vector space over $\F$ so that linear
transformations act on the left.
\m

Fix sequences
${\mathbf r} = \{r_n\}_{n \geqq 1}$\,, ${\mathbf s} = \{s_n\}_{n \geqq 1}$\,,
plus two new sequences ${\mathbf t'} = \{t'_n\}_{n \geqq 1}$ and
${\mathbf t''} = \{t''_n\}_{n \geqq 1}$\,, all of non--negative integers
with each $r_n + s_n > 0$.
and $d''_{n+1} = d''_n(r_{n+1} + s_{n+1}) + t''_{n+1}$ and denote
$d_{n+1} = d'_{n+1} + d''_{n+1}$\,.  Let $G_n$ be the real special unitary
group $SU(d'_n,d''_n;\F)$ over $\F$.  If $\F = \H$, or if $\F = \R$ and 
$d_n$ is odd, then $G_{n,\C}$ has no outer automorphism, and we denote 
$\delta = 1 \in Aut(G_n)$.  Otherwise (except when $\F = \R$ and $d_n = 8$)
$G_{n,\C}$ has outer automorphism group generated modulo inner 
automorphisms by $\delta_0 := 
Ad\left ( \begin{smallmatrix} -1 & 0 \\ 0 & +I \end{smallmatrix}\right )$
if $\F = \R$, by $\delta_0: g \mapsto {}^tg^{-1}$ if $\F = \C$,
and we choose $\delta \in \delta_0 Int(G_n)$ that preserves the standard 
positive root system.  Then we have
$\phi_{n+1,n}: SU(d'_n,d''_n;\F) \to SU(d'_{n+1},d''_{n+1};\F)$ 
given by
\begin{equation}\label{ddl-indefunitarygroups}
\phi_{n+1,n}(g) \ = \diag\{1, \dots , 1; g, \dots , g; \delta(g) , \dots , 
\delta(g) ; 1, \dots , 1\}
\end{equation}
with $t'_{n+1}$ entries $1$, $r_{n+1}$ blocks $g$, $s_{n+1}$ blocks 
$\delta(g)$, and finally $t''_{n+1}$ entries $1$, where all $s_n = 0$
in the case $\F = \H$.  Now (\ref{ddl-indefunitarygroups}) defines a
strict direct system $\{G_m,\phi_{n,m}\}$.  Let $d' = \lim d'_n$ and
$d'' = \lim d''_n$\,; both usually are $\infty$ but of course it can happen 
that one is finite, even zero.  In any case we have
\begin{equation} \label{ddl-indefunitarygroups-lim}
G = SU_{{\mathbf t'},{\mathbf r},{\mathbf s},{\mathbf t''}}(d',d'';\F) := 
\varinjlim \{G_m,\phi_{n,m}\} \text{ for the given } {\mathbf t'}, 
{\mathbf r}, {\mathbf s} \text{ and } {\mathbf t''}.
\end{equation}
Thus we have the groups 
$SO_{{\mathbf t'},{\mathbf r},{\mathbf s},{\mathbf t''}}(d',d'')$,
$SU_{{\mathbf t'},{\mathbf r},{\mathbf s},{\mathbf t''}}(d',d'')$ and
$Sp_{{\mathbf t'},{\mathbf r},{\mathbf s},{\mathbf t''}}(d',d'')$.
\m

The same process gives us real limit orthogonal groups
$O_{{\mathbf t'},{\mathbf r},{\mathbf s},{\mathbf t''}}(d',d'')$
and the complex limit unitary groups 
$U_{{\mathbf t'},{\mathbf r},{\mathbf s},{\mathbf t''}}(d',d'')$.  In the 
$O_{{\mathbf t'},{\mathbf r},{\mathbf s},{\mathbf t''}}(d',d'')$ 
case, (\ref{group-class1}) requires that each $d_n$ should be odd.
 $s_n = 0$.
\m

{\sc Symplectic groups.}
We consider limits of real and complex symplectic groups.  Fix sequences
${\mathbf r} = \{r_n\}_{n \geqq 1}$ and ${\mathbf t} = \{t_n\}_{n \geqq 1}$
with all $r_n > 0$.
Start with $d_0 > 0$ and recursively define $d_{n+1}=d_nr_{n+1} + t_{n+1}$\,.
Our convention is that $Sp(n;\F)$ is the automorphism group of $\F^{2n}$
with a nondegenerate antisymmetric bilinear form; that forces $\F$ to be
$\R$ or $\C$.  Let $G_n = Sp(d_n;\F)$, either the
real symplectic group $Sp(d_n;\R)$ or the complex symplectic group
$Sp(n;\C)$.  Then we have strict direct systems 
$\{G_m,\phi_{n,m}\}_{m \geqq n \geqq 0}$ with
\begin{equation} \label{ddl-symplecticgroups}
\phi_{n+1,n}: G_n \to G_{n+1} \text{ by }
\phi_{n+1,n}(g) \ = \diag\{g, \dots , g; 1, \dots , 1\}
\end{equation}
with $r_{n+1}$ blocks $g$ and with $2t_{n+1}$ entries $1$.  That gives us
\begin{equation} \label{ddl-sympleticgroups-lim}
G = Sp_{{\mathbf r},2{\mathbf t}}(\infty; \F) :=
\varinjlim \{G_m,\phi_{n,m}\} \text{ for the given } {\mathbf r}
\text{ and } {\mathbf t}.
\end{equation}

{\sc Complex orthogonal groups.}
Now consider the complex special orthogonal groups $G_n = SO(d_n;\C)$.
The formula (\ref{ddl-lineargroups-lim}) defines maps
$\phi_{n+1,n}: SO(d_n;\C) \to SO(d_{n+1};\C)$ so it defines a
strict direct system $\{G_m,\phi_{n,m}\}$.  Then we have
\begin{equation} \label{ddl-complexorthoggroups-lim}
G = SO_{{\mathbf r},{\mathbf s},{\mathbf t}}(\infty;\C) :=
\varinjlim \{G_m,\phi_{n,m}\} \text{ for the given } {\mathbf r}, {\mathbf s},
\text{ and } {\mathbf t}.
\end{equation}

The same process gives us  complex limit orthogonal groups 
$O_{{\mathbf r},{\mathbf s},{\mathbf t}}(\infty;\C)$; as before, 
here (\ref{group-class1}) requires that each $d_n$ be odd.
\m

{\sc The remaining classical series.}
There is one other series of real classical groups, the groups $SO^*(2n)$,
real form of $SO(2n;\C)$ with maximal compact subgroup $U(n)$.
The usual definition is 
$$
SO^*(2n) = \{g \in U(n,n) \mid g \text{ preserves }
(x,y) = {\sum}_1^n (x_a y_{n+a} + x_{n+a} y_a ) \text{ on } \C^{2n}\}.
$$
It will be more convenient for us to use the alternate formulation of
\cite[Section 8]{W3}, which is
\begin{equation} \label{sostar}
SO^*(2n) = \{g \in SL(n;\H) \mid b(gx,gy) = b(x,y) \text{ for all }
x,y \in \H^n\}
\end{equation}
where $b$ is the skew--hermitian form on $\H^n$ given by 
$b(x,y) = {\sum}_{a=1}^n \overline{x}_aiy_a$\,.  For then
(\ref{outer}) defines an outer automorphism $\delta$ of each
$SO^*(2n)$.  Now fix sequences ${\mathbf r} = \{r_n\}_{n \geqq 1}$,
${\mathbf s} = \{s_n\}_{n \geqq 1}$ and
${\mathbf t} = \{t_n\}_{n \geqq 1}$ of non--negative integers with
$r_0 > 0$ and all $r_n + s_n > 0$.  Start with  $d_0 > 0$ and recursively 
define $d_{n+1} = d_n(r_{n+1} + s_{n+1}) + t_{n+1}$\,.  Define
$G_n = SO^*(2d_n)$.  Then we have strict
direct systems $\{G_m,\phi_{n,m}\}_{n \geqq m \geqq 0}$ given by
\begin{equation}\label{ddl-sostargroups}
\phi_{n+1,n}: G_n \to G_{n+1} \text{ by }
\phi_{n+1,n}(g) \ = \diag\{g, \dots , g; \delta(g) , \dots , \delta(g) ;
1, \dots , 1\}
\end{equation}
with $r_{n+1}$ blocks $g$, with $s_{n+1}$ blocks $\delta(g)$, and with
$t_{n+1}$ entries $1$.  That gives us
\begin{equation} \label{ddl-sostargroups-lim}
G = SO^*_{{\mathbf r},{\mathbf s},{\mathbf t}}(\infty) :=
\varinjlim \{G_m,\phi_{n,m}\} \text{ for the given } {\mathbf r}, {\mathbf s}
\text{ and } {\mathbf t}.
\end{equation}

We refer to the direct limit groups (\ref{ddl-lineargroups-lim}),
(\ref{ddl-indefunitarygroups-lim}), (\ref{ddl-sympleticgroups-lim}),
(\ref{ddl-complexorthoggroups-lim}),
and (\ref{ddl-sostargroups-lim}) as {\em diagonal embedding direct limit
groups} and to the associated direct systems as {\em diagonal embedding direct
systems}.  In the unitary symplectic case of (\ref{ddl-indefunitarygroups-lim})
we made the convention that we have the sequence ${\mathbf s}$ but each 
$s_n = 0$.
We say that a diagonal embedding direct limit group and the  associated
diagonal embedding direct system are of {\em classical type} if
$r_n + s_n = 1$ for all $n$ sufficiently large.
\m

Now we collect some basic properties of diagonal embedding direct limit groups.

\begin{proposition}\label{diag-lim-exists}
Let $G = \varinjlim \{G_m,\phi_{n,m}\}_{n \geqq m \geqq 0}$ be a diagonal
embedding direct limit group.  Then conditions {\em (\ref{group-class1}) and
(\ref{group-class2})} always hold, and if the $G_n$ are not 
{\rm (}special{\rm )} unitary groups over $\R$, $\C$ or $\H$ then
{\rm (\ref{coherent-iwasawa})} holds, so $G$ has principal series 
representations.  If
the $G_n$ are {\rm (}special{\rm )} unitary groups over $\R$, 
$\C$ or $\H$, then {\rm (\ref{coherent-iwasawa})} holds  for a cofinal
subsystem {\rm (}which of course yields the same limit group $G${\rm )}
of $\{G_m,\phi_{n,m}\}$.  In any case, if $\{G_m,\phi_{n,m}\}$ is weakly 
parabolic then it is of classical type.
\end{proposition}

\noindent {\bf Proof.}  Conditions (\ref{group-class1}) and 
(\ref{group-class2}) are clear: the $G_n$ are semisimple Lie groups,
connected except possible for the case of orthogonal groups where we 
have explicitly ensured (\ref{group-class1}).  Now we look at 
(\ref{coherent-iwasawa}).
\m

We first consider the special linear groups $G_n = SL(d_n;\F)$.
Fix a basis $\cB$ of $\F^{d_n}$.  Relative to
$\cB$, $A_n$ will consist of the diagonal real matrices in
$G_n$ that have all entries $> 0$, $M_n$ will consist of the diagonal matrices 
in $G_n$ that have all entries of absolute value $1$, and $N_n$ will consist
of all lower triangular matrices in $G_n$ that have all diagonal entries
$= 1$.  It is immediate that $\phi_{n+1,n}$ maps $A_n$ into
$A_{n+1}$\,, maps $M_n$ into $M_{n+1}$\,, and maps $N_n$ into $N_{n+1}$\,.
That is (\ref{coherent-iwasawa}).
\m

Now consider the symplectic groups $G_n = Sp(d_n;\F)$.  The standard
basis $\{e_i\}$ of $\F^{2d_n}$, in which the antisymmetric bilinear form
$b_n$ that defines $G_n$ has matrix $\left ( \begin{smallmatrix} O & I \\ 
-I & 0 \end{smallmatrix} \right )$, specifies a new basis 
$\cB = \{v_1 , \ldots v_{d_n};
v'_1 , \ldots , v'_{d_n}\}$ by $v_i = \frac{1}{\sqrt{2}}(e_i + e_{d_n+i})$
and $v'_i = \frac{1}{\sqrt{2}}(e_i - e_{d_n+i})$.  Relative to $\cB$, the 
group $A_n$ will consist of the diagonal real matrices in $G_n$ with all
entries $>0$, in other words all $\diag\{a_1,\ldots , a_{d_n}, a_1^{-1},
\ldots , a_{d_n}^{-1}\}$ with each $a_i > 0$.  Then, as above,
$M_n$ will consist of the diagonal matrices
in $G_n$ that have all entries of absolute value $1$, and $N_n$ will consist
of all lower triangular matrices in $G_n$ that have all diagonal entries
$= 1$, so $\phi_{n+1,n}$ maps $A_n$ into
$A_{n+1}$\,, maps $M_n$ into $M_{n+1}$\,, and maps $N_n$ into $N_{n+1}$\,.
That is (\ref{coherent-iwasawa}).
\m

Next consider the complex special orthogonal groups $G_n = SO(d_n;\C)$.
Let $m_n = [d_n/2]$, let $\{e_i\}$ be a basis of $\C^{d_n}$ in which the
symmetric bilinear form $b_n$ that defines $G_n$ has matrix $I$.
Define $v_i = \frac{1}{\sqrt{2}}(e_i + e_{m_n+i}{\mathbf i})$ and
$v'_i = \frac{1}{\sqrt{2}}(e_i - e_{m_n+i}{\mathbf i})$ for 
$1 \leqq i \leqq m_n$\,.
Consider the basis $\cB$ of $\C^{d_n}$ given by
$\{v_1 , \ldots v_{m_n}; v'_1 , \ldots , v'_{m_n}\}$ if $d_n$ is even 
(hence $= 2m_n$),
by $\{v_1 , \ldots v_{m_n}; v'_1 , \ldots , v'_{m_n}; e_{d_n}\}$ if $d_n$ 
is odd (hence $= 2m_n+1$).  Relative to $\cB$, if $d_n$ is even then
$A_n$ will consist of all matrices $\diag\{a_1,\ldots , a_{m_n}, a_1^{-1},
\ldots , a_{m_n}^{-1}\}$ with each $a_i > 0$, and if $d_n$ is odd it will 
consist of all 
$\diag\{a_1,\ldots , a_{m_n}, a_1^{-1}, \ldots , a_{m_n}^{-1},1\}$.  Then
$M_n$ will consist of the diagonal matrices
in $G_n$ that have all entries of absolute value $1$, and $N_n$ will consist
of all lower triangular matrices in $G_n$ that have all diagonal entries
$= 1$, so $\phi_{n+1,n}$ maps $A_n$ into
$A_{n+1}$\,, maps $M_n$ into $M_{n+1}$\,, and maps $N_n$ into $N_{n+1}$\,.
That is (\ref{coherent-iwasawa}).
\m

We now consider the groups $G_n = SO^*(2d_n)$, essentially as above.
Let $U_n$ be the underlying right vector space over $\H$ on which
$G_n$ acts.  Let $\{e_i\}$ be a basis of $U_n = \H^{d_n}$ in which the
skew--hermitian form $b_n$ that defines $G_n$ is given by
$b_n(z,w) = \sum_{1 \leqq i \leqq d_n} \overline{w_i}{\mathbf i}z_i$\,.
Define $v_i = \frac{1}{\sqrt{2}}(e_{2i-1} + e_{2i}{\mathbf j})$ and
$v'_i = \frac{1}{\sqrt{2}}(e_{2i-1}{\mathbf i} + e_{2i}{\mathbf k})$ 
for $1 \leqq i \leqq m_n$ where $m_n = [d_n/2]$.  Let
$\cB = \{v_1, \ldots , v_{m_n}; v'_1, \ldots v'_{m_n}\}$ if $d_n$ is even, i.e.
$d_n = 2m_n$\,, and let $\cB = 
\{v_1, \ldots , v_{m_n}; v'_1, \ldots v'_{m_n};e_{d_n}\}$ if $d_n$ is odd,
i.e. $d_n = 2m_n+1$\,.  Then $V_n = \text{ Span}\{v_i\}$ and
$V'_n = \text{Span}\{v'_i\}$ are maximal totally $b_n$--isotropic subspaces
of $U_n$\,, paired by $b_n(v_i,v'_j) = \delta_{ij}$.  If $d_n$ is even
then $U_n = V_n + V'_n$\,, and if $d_n$ is odd then $U_n = V_n + V'_n + W_n$
where $W_n = \text{ Span}\{e_{d_n}\} = (V_n + V'_n)^\perp$ relative to $b_n$\,.
In the basis $\cB$ the groups $A_n$, $M_n$ and $N_n$ are given as in
the case of the complex special orthogonal groups, so $\phi_{n+1,n}$ maps 
$A_n$ into $A_{n+1}$\,, maps $M_n$ into $M_{n+1}$\,, and maps $N_n$ into 
$N_{n+1}$\,.  That gives us (\ref{coherent-iwasawa}).
\m

Finally we come to the case $G_n = SU(d'_n,d''_n;\F)$ of the real orthogonal, 
complex unitary and unitary symplectic (quaternion unitary) groups.
Let $U_n$ be the underlying right vector space, over $\F$ on which 
$G_n$ acts.  Then $G_n$ is essentially the group of automorphisms
of $(U_n,b_n)$ where $b_n$ is the nondegenerate $\F$--hermitian form on
$U_n$ that defines $G_n$\,.  Let $\{e_i\}$ be a basis of $U_n$ in which
$b_n$ has matrix $\left ( \begin{smallmatrix} I & 0 \\ 0 & -I 
\end{smallmatrix} \right )$.  Let $m_n = \min(d'_n,d''_n)$, the real
rank $\dim \ga_n$ of $G_n$\,.  Define
$v_i = \frac{1}{\sqrt{2}}(e_i + e_{m_n+i})$ and
$v'_i = \frac{1}{\sqrt{2}}(e_i - e_{m_n+i})$ for $1 \leqq i \leqq m_n$\,.
Let $r_n = d_n - 2m_n$\,, and let $\{w_1 , \ldots , w_{r_n}\}$ denote
the ordered set of those $e_i$ not involved in the $v_j$\,.  
Then we have the basis $\cB = \{v'_1, \ldots , v'_{m_n}; w_1 , \ldots w_{r_n};
v_{m_n}, \dots , v_1\}$ of $U_n$\,.  On the subspace level, 
$V_n = \text{ Span}\{v_i\}$ and
$V'_n = \text{Span}\{v'_i\}$ are maximal totally $b_n$--isotropic subspaces
of $U_n$\,, paired by $b_n(v_i,v'_j) = \delta_{ij}$, and 
$U_n = V'_n + W_n + V_n$ where $W_n = \text{ Span}\{w_i\} = 
(V_n + V'_n)^\perp$ relative to $b_n$\,.  Note that $W_n$ is zero if
$d'_n = d''_n$\,, positive
definite if $d'_n > d''_n$\,, negative definite if $d'_n < d''_n$\,.
\s

We choose $A_n$ to consist of all linear transformations of $U_n$ with
matrix, relative to $\cB$, of the form $\diag\{a_1,\ldots , a_{m_n}; 
1, \ldots 1; a_{m_n}^{-1}, \ldots , a_1^{-1}\}$ with $a_i$ all real and 
positive.  Then $M_n$ consists of all linear transformations 
$m \in G_n$ such that
$$
\begin{aligned}
&m(v'_i) = m_iv'_i\ , \ m(W_n) = W_n \ \text{ and }  m(v_i) = m_iv_i\  \ \ 
\text{where the } m_i \in \F \text{ with } |m_i| = 1 
\end{aligned}
$$ 
Thus $\phi_{n+1,n}(A_n) \subset A_{n+1}$\,.
\m

The description of $N_n$ is a little more complicated.  Let 
$\cV_n = \{V_{n,1}, \ldots , V_{n,m_n}\}$ be the maximal isotropic flag in
$V_n$ given by $V_{n,j} = \text{ Span}\{v_1, \ldots , v_j\}$.
In almost every case we may take the minimal parabolic subgroup 
$P_n$ of $G_n$ to be the $G_n$--stabilizer of $\cV_n$\,.  That done, 
let $P_{n,j}$ denote the maximal real parabolic subgroup of $G_n$ that
is the stabilizer of $V_{n,j}$.  Then the nilradicals of these parabolics
satisfy $\gn_n = \sum_{1\leqq j\leqq m_n} \gn_{n,j}$\,.  The point of this
is that we know the $\gn_{n,j}$ in a convenient form.
Let $X_{n,j} = \text{ Span}\{v_{j+1} , \ldots , v_{m_n}\}$ so that 
$V_n = V_{n,j} \oplus X_{n,j}$\,.  Define 
$V'_{n,j} = \text{ Span}\{v'_1, \ldots , v'_j\}$ and
$X'_{n,j} = \text{ Span}\{v'_{j+1} , \ldots , v'_{m_n}\}$ so that
$V'_n = V'_{n,j} \oplus X'_{n,j}$\,.  Denote 
$W_{n,j} =  X_{n,j} \oplus W_n \oplus X'_{n,j}$ so that 
$U_n = V_{n,j} + W_{n,j} + V'_{n,j}$\,.  According to \cite[Lemma 3.4]{W2},
the nilradical $\gn_{n,j}$ of $\gp_{n,j}$ is the sum of its two subspaces
$$
\begin{aligned}
&\gp_{n,j}^2 = \{\xi \in \gg_n \mid \xi(V'_{n,j}) 
	\subset V_{n,j}\,, \xi(W_{n,j}) = 0\,, \xi (V_{n,j}) = 0\}
	\text{ and } \\
&\gp_{n,j}^1 = \{\xi \in \gg_n \mid \xi(V'_{n,j}) \subset W_{n,j}\,, 
	\xi (W_{n,j}) \subset V_{n,j}\,, \xi (V_{n,j}) = 0\}
\end{aligned}
$$
while the reductive component consists of those $\xi$ in $\gg_n$ that
stabilize each of $V'_{n,j}$\,, $W_{n,j}$ and $V_{n,j}$\,.  Thus, relative
to the basis $\cB$, the elements of $\gn_{n,j}$ have block form
$\left ( \begin{smallmatrix} 0 & 0 & 0 \\ * & 0 & 0 \\ * & * & 0
\end{smallmatrix} \right )$ along $U_n = V'_{n,j} + W_{n,j} + V_{n,j}$\,.
Summing over $j$, the elements of $\gn_n$ are precisely those elements
of $\gp_n$ whose matrix relative to $\cB$ has block form
$\left ( \begin{smallmatrix} \ell & 0 & 0 \\ * & 0 & 0 \\ * & * & \ell'
\end{smallmatrix} \right )$ along $U_n = V'_n + W_n + V_n$\,, where 
$\ell$ and $\ell'$ are lower triangular with zeroes on their diagonals.
In the case $G_n = SO(d'_n,d''_n)$ one must be a bit more careful and
take some orientation into account, as in \cite{WZ}, but the result is
the same.  Thus $\phi_{n+1,n}(N_n) \subset N_{n+1}$\,.
\m

It certainly cannot be automatic that $\phi_{n+1,n}(M_n) \subset M_{n+1}$\,.
Define a difference $\mu_n = \max(d'_n,d''_n) - \min(d'_n,d''_n)$.  Then
$\mu_n = \dim W_n$\,, so $M^0_n \cong SU(\mu_n;\F)$ and 
$M_{n+1}^0 \cong SU(\mu_{n+1};\F)$, and thus 
$\phi_{n+1,n}(M_n) \not \subset M_{n+1}$ whenever $\mu_n > \mu_{n+1}$\,.
\m

Now we pin this down.
The map $\phi_{n+1,n}: G_n \to G_{n+1}$ is implemented by a unitary 
injection $k_n : (U_n,b_n) \hookrightarrow (U_{n+1},\pm b_{n+1})$\,.  
We have set things up so  that, possibly after interchanging the
$v_i$ and the $v'_i$ in $U_{n+1}$\,, $k_n(V_{n,j}) \subset V_{n+1,j}$ and 
$k_n(V'_{n,j}) \subset V'_{n+1,j}$ for $1 \leqq j \leqq m_n$\,.
We used that to prove that
$\phi_{n+1,n}$ maps $A_n$ into $A_{n+1}$ and $N_n$ into $N_{n+1}$\,.
But $\phi_{n+1,n}(M_n) \subset M_{n+1}$ if and only if we can make the
choices of $V_{n+1}$ and $V'_{n+1}$ so that $k_n(W_n) \hookrightarrow
W_{n+1}$\,.  That is possible if and only if $\mu_n \leqq \mu_{n+1}$\,.
\m

If $\mu_n \leqq \mu_{n+1}$ for infinitely many indices $n$, then we have a
cofinal subsystem of $\{G_m,\phi_{n,m}\}$ in which $\mu_n \leqq \mu_{n+1}$
for all $n$, and thus $\phi_{n+1,n}(M_n) \subset M_{n+1}$\,.  If
$\mu_n \leqq \mu_{n+1}$ for only finitely many indices $n$, then we have an
index $n_0$ such that $\mu_n > \mu_{n+1} \geqq 0$ for all $n \geqq n_0$.
That is impossible.  Thus  (\ref{coherent-iwasawa}) is always valid for
a cofinal subsystem of $\{G_m,\phi_{n,m}\}$.
\m

Suppose that $\{G_m,\phi_{n,m}\}$ is weakly parabolic.  View the
$\phi_{n+1,n}$ as inclusions $G_n \hookrightarrow G_{n+1}$\,.
Then $\ga_{n+1} = \ga_n \oplus \ga_{n+1,n}$ as in Proposition
\ref{rho-restriction}, and $(\ga_n + \gm_n) \oplus \ga_{n+1,n}$ 
is the centralizer of $\ga_{n+1,n}$ in $\gg_{n+1}$\,.  In
particular, $\Sigma(\gg_n \oplus \ga_{n+1,n}, \ga_{n+1})
\subset \Sigma(\gg_{n+1}, \ga_{n+1})$.  Thus, if $\gamma_n \in
\Sigma(\gg_n , \ga_n)$ there is a unique $\gamma_{n+1} \in
\Sigma(\gg_{n+1}, \ga_{n+1})$ such that $\gamma_{n+1}|_{\ga_n}
= \gamma_n$\,.  But if $r_{n+1} + s_{n+1} \geqq 2$ then at least two
distinct elements of $\Sigma(\gg_{n+1}, \ga_{n+1})$ restrict to
$\gamma_n$\,.  Thus $r_{n+1} + s_{n+1} = 1$.
\hfill $\square$

\begin{proposition}\label{class-cond}
Let $G = \varinjlim \{G_m,\phi_{n,m}\}_{n \geqq m \geqq 0}$ be a diagonal
embedding direct limit group.  Then the following conditions are
equivalent.

{\rm 1.} $G = \varinjlim \{G_m,\phi_{n,m}\}$ is of classical type, in other 
words $r_n + s_n = 1$ for $n \geqq n_0$\,.

{\rm 2.} The root system $\Sigma(\gg,\ga) = \varprojlim \Sigma(\gg_n, \ga_n)$
is countable.

{\rm 3.} $\Sigma(\gg,\ga) = {\bigcup}_{n \geqq 0}\Sigma(\gg_n, \ga_n)$.

{\rm 4.} $\gg$ is restricted--root--reductive in the sense that
$\gg = (\gm + \ga) + {\sum}_{\gamma \in \Sigma(\gg,\ga)}\, \gg^\gamma$\,.
\end{proposition}

\noindent {\bf Proof.}  Let $G$ be of classical type.  View the
$\phi_n$ as inclusions $G_n \hookrightarrow G$ as inclusions.  Let $\Psi_n$
denote the simple system of $\Sigma(\gg_n,\ga_n)^+$.  If $n \geqq m \geqq n_0$
then Proposition \ref{diag-lim-exists} shows how $\Psi_m \subset \Psi_n$
when we extend the elements of $\Psi_m$ by zero on $\ga_{n,m}$\,.
Thus $\Sigma(\gg,\ga) = {\bigcup}_{n \geqq 0}\Sigma(\gg_n, \ga_n)$, which is
countable, and if $\gamma \in \Sigma(\gg,\ga)$ then 
$\phi_n(\gg_n^{\phi_n^*(\gamma)})$ is a well defined subspace of the
root space for $\gamma$.  We have just seen that (1) implies (2), (3), and (4).
On the other hand, (4) implies (3), and (3) implies (2), at a glance.  Thus
we need only prove that (2) implies (1).
\m

Suppose that $G = \varinjlim \{G_m,\phi_{n,m}\}$ is not of classical type.  
Then we
can pass to a cofinal subsystem in which every $r_n + s_n \geqq 2$.  That
done, every root $\gamma \in \Sigma(\gg_n,\ga_n)^+$ is the restriction of
at least $2$ roots in $\Sigma(\gg_{n+1},\ga_{n+1})^+$, thus is the
restriction of at least $2^{\aleph_0}$ roots in $\Sigma(\gg,\ga)$.  In
particular $\Sigma(\gg,\ga)$ is not countable.  Thus (2) implies (1).
\hfill $\square$
\m

Recall the notion of Satake diagram.  We use the Cartan subalgebra 
$\gh_n = \gt_n + \ga_n$ of $\gg_n$ and the positive root system
$\Sigma(\gg_{n,\C},\gh_{n,\C})^+$ as in (\ref{coherent-root-order}).  Write
$\Psi(\gg_{n,\C},\gh_{n,\C})$ for the corresponding simple 
$\gh_{n,\C}$--root system,
and write $\Psi(\gg_n, \ga_n)$ for the simple $\ga_n$--root system
corresponding to $\Sigma(\gg_n,\ga_n)^+$.
Every $\psi \in \Psi(\gg_n, \ga_n)$ is of the form $\widetilde{\psi}|_{\ga_n}$
for some $\widetilde{\psi} \in \Psi(\gg_{n,\C},\gh_{n,\C})$.  More or less
conversely, if $\widetilde{\psi} \in \Psi(\gg_{n,\C},\gh_{n,\C})$ then
either $\widetilde{\psi}|_{\ga_n} = 0$ or $\widetilde{\psi}|_{\ga_n} \in
\Psi(\gg_n, \ga_n)$.  The Satake
diagram describes the restriction process.  Start with the Dynkin diagram
$\cD_n$ of $\gg_{n,\C}$ whose vertices are the elements of 
$\Psi(\gg_{n,\C},\gh_{n,\C})$.  If there are two root lengths this is
indicated by arrows rather than darkening the vertices for short roots.
Now darken those $\widetilde{\psi} \in \Psi(\gg_{n,\C},\gh_{n,\C})$ such that
$\widetilde{\psi}|_{\ga_n} = 0$.  It can happen that two (but never
more than two) distinct elements 
$\widetilde{\psi}, \widetilde{\psi}' \in \Psi(\gg_{n,\C},\gh_{n,\C})$
have the same $\ga_n$--restriction.  In that case, join them by an arrow.  
The result is the Satake diagram of $\gg_n$\,.  The white vertices  and
vertex pairs corresponding to simple $\ga_n$--roots of $\gg_n$.  The 
black vertices correspond to simple $\gt_{n\C}$--roots of $\gm_{n,\C}$\,.
See \cite{W4}, pp. 90--93 or \cite{A}, pp. 32--33, for Araki's list of 
Satake diagrams.
\m

We use the Satake diagrams to see just which
$G = \varinjlim \{G_m,\phi_{n,m}\}$
of classical type are weakly parabolic.  
The description (\ref{real-parabolics}) of real parabolic subalgebras
gives us
\begin{lemma}\label{satake}
The semisimple components of real parabolic subalgebras of $\gg_{n+1}$
are characterized up to $Int(G_{n+1})$--conjugacy by their Satake
diagrams, and those Satake diagrams are obtained from the Satake
diagram of $\gg_{n+1}$ by deleting
{\rm (i)} an arbitrary set of white vertices, and then {\rm (ii)} all
white vertices joined by arrows $($meaning the same restriction to 
$\ga_{n+1})$ to vertices deleted in {\rm (i)}.
\end{lemma}

The black vertices (restriction $0$ to $\ga_{n+1}$) remain because they
represent the simple roots of $\gm_{n+1}$\,, which is contained in every
real parabolic subalgebra that contains $\ga_{n+1}$\,.
\m

Let $G = \varinjlim \{G_m,\phi_{n,m}\}$ be a diagonal embedding direct limit 
group of classical type.  From Araki's list of Satake diagrams
one sees that the possible inclusions 
$\phi_{n+1,n}:\gg'_n \to \gg_{n+1}$ of weakly parabolic type are given,
modulo $\gm_{n+1}$, by
\m
\addtocounter{equation}{1}

\noindent (\theequation a)
$SL(d_n;\F) \hookrightarrow SL(d_{n+1},\F)$ by $\phi_{n+1,n}(g) =
	\diag\{g\text{ or } \delta(g),1,\ldots 1\}, d_{n+1} > d_n$\,, 
\m

\noindent (\theequation b)
$SO(d'_n,d''_n) \hookrightarrow SO(d'_n+u_n,d''_n+u_n)$
	by $\phi_{n+1,n}(g) =
        \diag\{g\text{ or } \delta(g),1,\ldots 1\}, u_n > 0$,
\m

\noindent (\theequation c)
$SO(d_n;\C) \hookrightarrow SO(d_n + 2u_n;\C)$ by $\phi_{n+1,n}(g) =
	\diag\{g\text{ or } \delta(g),1,\ldots 1\}, u_n > 0$,
\m

\noindent (\theequation d)
$SU(d'_n,d''_n) \hookrightarrow SU(d'_n+u_n,d''_n+u_n)$ by $\phi_{n+1,n}(g) =
        \diag\{g\text{ or } \delta(g),1,\ldots 1\}, u_n > 0$,
\m
 
\noindent (\theequation e)
$Sp(d'_n,d''_n) \hookrightarrow Sp(d'_n+u_n,d''_n+u_n)$ by $\phi_{n+1,n}(g) =
        \diag\{g,1,\ldots 1\}, u_n > 0$,
\m

\noindent (\theequation f)
$Sp(d_n;\F) \hookrightarrow Sp(d_{n+1};\F)$ by $\phi_{n+1,n}(g) =
	\diag\{g,1,\ldots 1\}, d_{n+1} > d_n$ and $\F = \R$ or $\C$,
\m

\noindent (\theequation g)
$SO^*(2d_n) \hookrightarrow SO^*(2d_n+4u_n)$  by $\phi_{n+1,n}(g) =
        \diag\{g\text{ or } \delta(g),1,\ldots 1\}, u_n > 0$.
\m

In order to pin things down we make use of the fact that 
$G = \varinjlim \{G_m,\phi_{n,m}\}$ is determined by any cofinal subsequence
of indices.  Denote
$$
{\mathbf 0} = \{0,0,0,\ldots\}, \ {\mathbf 1} = \{1,1,1,\ldots\}
\text{ and } {\mathbf 2} = \{2,2,2,\ldots\}.
$$
Consider, for example, the case of (\theequation a).  Suppose first that 
there are only finitely many indices $n$ for which $\phi_{n+1,n}(g) = 
\diag\{\delta(g),1, \ldots , 1\}$.  Pass to the subsequence starting just
after the last $\phi_{n+1,n}$ that involves $\delta$.  That done,
we interpolate and arrive at the same limit with
each $\phi_{n+1,n}(g) = \left ( \begin{smallmatrix} g & 0 \\ 0 & 1
\end{smallmatrix} \right )$.  Now suppose that there are infinitely
many indices $n$ for which $\phi_{n+1,n}(g) =
\diag\{\delta(g),1, \ldots , 1\}$.  Pass to the cofinal subsequence
obtained by deleting the $G_n$ for which $\phi_{n+1,n}(g) =
\diag\{g,1,\ldots ,1\}$, so now every $\phi_{n+1,n}(g)$ is of the form
$g \mapsto \diag\{\delta(g),1, \ldots , 1\}$.  If $t_{n+1} > 1$ for an 
infinite number of $t_{n+1}$ then, recursively, we take the smallest index $n$
for which $t_{n+1} > 1$, insert $t_{n+1} - 1$ steps $g \mapsto 
\left ( \begin{smallmatrix} \delta(g) & 0 \\ 0 & 1 \end{smallmatrix} \right )$
between $G_n$ and $G_{n+1}$\,, and proceed to insert steps $g \mapsto
\left ( \begin{smallmatrix} \delta(g) & 0 \\ 0 & 1 \end{smallmatrix} \right )$
at the next $t_{n+1} - 1$ possible places.  Then we arrive at the same
limit with each $\phi_{n+1,n}(g) = \left ( \begin{smallmatrix} \delta(g) & 0 
\\ 0 & 1 \end{smallmatrix} \right )$.  If $t_{n+1} > 1$ for
only finitely many $n$ we just pass to the subsequence starting just
after the last $\phi_{n+1,n}$ involving a $t_{n+1}$ that is $> 1$.  Thus
$G = SL_{{\mathbf 1},{\mathbf 0},{\mathbf 1}}(\infty,\F) = 
\varinjlim SL(n+1;\F)$ in the first case, with $\phi_{n+1,n}(g) =
\left ( \begin{smallmatrix} g & 0 \\ 0 & 1 \end{smallmatrix} \right )$, and 
$G = SL_{{\mathbf 0},{\mathbf 1},{\mathbf 1}}(\infty,\F) = 
\varinjlim SL(n+1;\F)$ in the second case, with $\phi_{n+1,n}(g) =
\left ( \begin{smallmatrix} \delta(g) & 0 \\ 0 & 1 \end{smallmatrix} \right )$.
\m

Similar considerations hold in the other six cases of (\theequation).  The
final result is

\begin{proposition}\label{classical-parabolic}
The weakly parabolic diagonal embedding direct limit groups
$G = \{G_m,\phi_{n,m}\}$ of classical type, with $G_m$ noncompact and simple 
for $m$ large, are given, up to isomorphism, by one of the following.
\m
\addtocounter{equation}{1}

\noindent {\rm (\theequation a)} 
$SL_{{\mathbf 1},{\mathbf 0},{\mathbf 1}}(\infty;\F)$ with $g \mapsto 
\left ( \begin{smallmatrix} g & 0 \\ 0 & 1 \end{smallmatrix} \right )$ or
$SL_{{\mathbf 0},{\mathbf 1},{\mathbf 1}}(\infty,\F)$  with $g \mapsto
\left ( \begin{smallmatrix} \delta(g) & 0 \\ 0 & 1 \end{smallmatrix} \right )$;
\hfill\newline\phantom{XXXX} here we may take $G_m$ to be $SL(m+1;\F)$.
\m

\noindent {\rm (\theequation b)}
$SO_{{\mathbf 1},{\mathbf 1},{\mathbf 0},{\mathbf 1}}(\infty,\infty)$ with
$g \mapsto \left ( \begin{smallmatrix} 1 & 0 & 0 \\ 0 & g & 0 \\ 0 & 0 & 1 
\end{smallmatrix} \right )$ or 
$SO_{{\mathbf 1},{\mathbf 0},{\mathbf 1},{\mathbf 1}}(\infty,\infty)$ with
$g \mapsto \left ( \begin{smallmatrix} 1 & 0 & 0 \\ 0 & \delta(g)  & 0 \\ 
0 & 0 & 1 \end{smallmatrix} \right )$;
\hfill\newline\phantom{XXXX} here we may take $G_m$ to be an 
$SO(d'_1+m,d''_1+m)$ where $d'_1, d''_1 \geqq 1$.
\m
 
\noindent {\rm (\theequation c)}
$SO_{{\mathbf 1},{\mathbf 0},{\mathbf 2}}(\infty;\C)$ with $g \mapsto
\left ( \begin{smallmatrix} g & 0 & 0 \\ 0 & 1 & 0 \\ 0 & 0 & 1
\end{smallmatrix} \right )$ or
$SO_{{\mathbf 0},{\mathbf 1},{\mathbf 2}}(\infty;\C)$ with $g \mapsto
\left ( \begin{smallmatrix} \delta(g) & 0 & 0 \\ 0 & 1 & 0 \\ 0 & 0 & 1
\end{smallmatrix} \right )$;
\hfill\newline\phantom{XXXX} here we may take $G_m$ to be
$SO(2m+1;\C)$ $($type B$)$ or $SO(2m;\C)$ $($type D$)$.
\m
 
\noindent {\rm (\theequation d)}
$SU_{{\mathbf 1},{\mathbf 1},{\mathbf 0},{\mathbf 1}}(\infty;\infty)$ 
with $g \mapsto
\left ( \begin{smallmatrix} 1 & 0 & 0 \\ 0 & g & 0 \\ 0 & 0 & 1
\end{smallmatrix} \right )$ or
$SU_{{\mathbf 1},{\mathbf 0},{\mathbf 1},{\mathbf 1}}(\infty;\infty)$ 
with $g \mapsto
\left ( \begin{smallmatrix} 1 & 0 & 0 \\ 0 & \delta(g) & 0 \\ 0 & 0 & 1
\end{smallmatrix} \right )$;
\hfill\newline\phantom{XXXX} here we may take $G_m$ to be an 
$SU(d'_1+m,d''_1+m)$ where $d'_1, d''_1 \geqq 1$.
\m
 
\noindent {\rm (\theequation e)}
$Sp_{{\mathbf 1},{\mathbf 1},{\mathbf 0},{\mathbf 1}}(\infty;\infty)$ with $g \mapsto
\left ( \begin{smallmatrix} 1 & 0 & 0 \\ 0 & g & 0 \\ 0 & 0 & 1
\end{smallmatrix} \right )$; 
\hfill\newline\phantom{XXXX} here we may take $G_m$ to be an 
$Sp(d'_1+m,d''_1+m)$ where $d'_1, d''_1 \geqq 1$.
\m
 
\noindent {\rm (\theequation f)}
$Sp_{{\mathbf 1},{\mathbf 2}}(\infty;\F)$ with $g \mapsto 
\left ( \begin{smallmatrix} g & 0 & 0 \\ 0 & 1 & 0 \\ 0 & 0 & 1
\end{smallmatrix} \right )$ and $\F = \R$ or $\C$; 
\hfill\newline\phantom{XXXX} here we may take $G_m$ to be $Sp(m;\F)$.
\m
 
\noindent {\rm (\theequation g)}
$SO^*_{{\mathbf 1},{\mathbf 0},{\mathbf 1}}(\infty)$ with $g \mapsto
\left ( \begin{smallmatrix} g & 0 \\ 0 & 1 \end{smallmatrix} \right )$ or
$SO^*_{{\mathbf 0},{\mathbf 1},{\mathbf 1}}(\infty)$  with $g \mapsto
\left ( \begin{smallmatrix} \delta(g) & 0 \\ 0 & 1 \end{smallmatrix} \right )$,
$($quaternionic matrices$)$;
\hfill\newline\phantom{XXXX} here we may take $G_m$ to be $SO^*(2m)$.
\end{proposition}

\section{The Other Tempered Series}
\label{sec10}
\setcounter{equation}{0}

The finite dimensional real reductive Lie groups $G$ that satisfy
(\ref{group-class1}) and (\ref{group-class2}) have a series of unitary
representations for each conjugacy class of Cartan subgroups.  Those
are the ``tempered'' representations, the ones that occur in the
decomposition of $L_2(G)$ under the left translation action of $G$.
The principal series is the tempered series corresponding to the
conjugacy class of a maximally noncompact Cartan subgroup, but of course
there are others.  If $G$ has a Cartan subgroup with compact image 
under the adjoint representation, then the corresponding series is 
the discrete series.  In general these series are constructed by 
combining the ideas of the discrete series and the principal series.
See \cite{H0}, \cite{H1}, \cite{H2} and \cite{H3} for the case where
$G$ is Harish--Chandra class, \cite{W1}, \cite{HW1} and \cite{HW2} for
the general case.  We now recall a few relevant facts from these papers
in order to indicate the corresponding extension of our principal series
results.
\m

Fix a Cartan involution $\theta$ of $G$ and let $K$ denote its fixed point
set.  As usual, $\gg = \gk + \gs$, decomposition into $(\pm 1)$--eigenspaces
of $\theta$.  Every $G^0$--conjugacy class of
Cartan subgroup contains a $\theta$--stable Cartan.  Fix a 
$\theta$--stable Cartan subgroup $H$ of $G$.  Then $\gh = \gt + \ga$ where
$\gt = \gh \cap \gk$ and $\ga = \gh \cap \gs$.  Here $H = T \times A$ where
$T = H \cap K$ and $A = \exp(\ga)$.  Earlier we had only considered the
case where $\ga$ is maximal abelian in $\gs$; here the situation is more
general.  The centralizer of $A$ in $G$ has form $M \times A$ where
$\theta(M) = M$.  In our earlier discussions $M$ was compact modulo
$Z_G(G^0)$ ({\em relatively compact}), but here the situation is more 
general.  In any case, $M$
satisfies (\ref{group-class1}) and (\ref{group-class2}), and $T$ is
relatively compact, so $M$ has relative discrete series representations.
In the principal series setting these will be all the irreducible
representations of $M$ and will necessarily be finite dimensional, but here
the situation is more general.  
\m

We have the $\ga$--root system
$\Sigma(\gg,\ga):= \{\alpha|_\ga \ \mid \ \alpha \in \Sigma(\gg,\gh)
\text{ and } \alpha|_\ga \ne 0\}$.  Fix a positive subsystem
$\Sigma(\gg,\ga)^+$ and define $\gn = \sum_{\beta \in \Sigma(\gg,\ga)^+} 
\gg^{-\beta}$.  Then $\gp := \gm + \ga + \gn$ is a particular kind of
(real) parabolic subalgebra of $\gg$, distinguished by the fact that
$\gt$ is a Cartan subalgebra of $\gm$.  Those are the {\em cuspidal
parabolic} subalgebras of $\gg$.  Let $N = \exp(\gn)$.  It is the analytic
subgroup of $G$ with Lie algebra $\gn$, and $P = MAN$ is the 
parabolic subgroup of $G$ with Lie algebra $\gp$.  Those are the
{\em cuspidal parabolic} subgroups of $G$.  
\m

We use the cuspidal parabolic subgroup $P = MAN$ to describe the
$H$--series representations of $G$.  The analog of Proposition 
\ref{m-structure} is the parameterization of the relative discrete
series of $M$.  Fix a positive $\gt_{_\C}$--root system
$\Sigma(\gm_{_\C},\gt_{_\C})^+$ on $\gm_{_\C}$\,.  Let $\nu \in i\gt^*$
such that $e^{\nu - \rho_{\gm,\gt}}$ is well defined on $T^0$ and
$\langle \nu , \alpha \rangle \ne 0$ for all 
$\alpha \in \Sigma(\gm_{_\C},\gt_{_\C})$.  Then $M^0$ has a unique
unitary equivalence class of relative discrete series representations,
$[\eta_\nu^0]$, with Harish--Chandra parameter $\nu$.  Here $\nu$ is the
infinitesimal character of $[\eta_\nu^0]$; if $\eta_\nu^0$ has a highest 
weight, that weight is $\nu - \rho_{\gm,\gt}$.  Set $M^\dagger = Z_M(M^0)M^0$.
Then the relative discrete series classes of $M^\dagger$ are the
$[\eta_{\chi,\nu}^\dagger] = [\chi \otimes \eta_\nu]$ where
$[\chi] \in \widehat{Z_M(M^0)}_\xi$ with $\xi = 
e^{\nu - \rho_{\gm,\gt}}|_{Z_{M^0}}$\,.  The relative discrete series classes
of $M$ are the $[\eta_{\chi,\nu}]$ where 
$\eta_{\chi,\nu} = \text{\rm Ind}_{M^\dagger}^M(\eta_{\chi,\nu}^\dagger)$.
Let $\sigma \in \ga_{_\C}^*$\,.  Then we have $[\eta_{\chi,\nu,\sigma}] 
\in \widehat{P}$ defined by 
$\eta_{\chi,\nu,\sigma}(man) = e^\sigma(a)\eta_{\chi,\nu}(m)$ for $m \in M,
a \in A \text{ and } n \in N$.  The corresponding $H$--series 
rep-resentation of $G$ is $\pi_{\chi,\nu,\sigma} = 
\text{\rm Ind}_P^G(\eta_{\chi,\nu,\sigma})$.  Its equivalence class does 
not depend on the choice of $\Sigma(\gg,\ga)^+$.  The $H$--series of $G$
consists of all such representations --- or, depending on context, the
unitary ones.  The principal series of $G$ is the 
case where $\ga$ is maximal, equivalently where $M$ is relatively compact.
The relative discrete series of $G$ is the case $\ga = 0$; it exists if
and only if $G$ has a relatively compact Cartan subgroup.
\m

Now we return to our strict direct system $\{G_i,\phi_{k,i}\}$
of reductive Lie groups.  Fix a Cartan subgroup $H = \varinjlim H_i$
of $G = \varinjlim G_i$\,.  We consider limit 
representations $\pi = \varprojlim \pi_i$ of $G$, where, for each
$i$, $\pi_i$ is an $H_i$--series representation of $G_i$\,.
\m

Here there are several problems.  First, we need the discrete series
analog of \cite{NRW3} in order to construct the $M$--component of any 
$\varprojlim \pi_i$\,.  That falls into two parts.  The first is to
realize the discrete series representations of the $M_i$ on some
appropriate cohomology spaces, such as spaces of $L_2$ harmonic forms.  
This is done, for example, in \cite{S1}, \cite{S2} and \cite{W1}.  The 
second is to
make sure that these representations all appear on cohomologies of the
same degree, and to line them up properly so that one can take limits.
This was done in \cite{N} for holomorphic discrete series; there the
cohomology degree is $0$, the alignment is done using the universal
enveloping algebra description of highest weight representations, and
the result is analysed by use of \cite{EHW}.  It was done in \cite{Hab}
for other discrete series of certain diagonal embedding direct limit
groups $Sp(p,\infty)$ and $SO(2p,\infty)$ of classical type using 
Zuckerman derived functor modules $A_q(\lambda)$ for the cohomologies.
We address these matters in some generality in \cite{W5}.
\m

Second, we need an analog of the considerations of Section \ref{sec8}.
This is not so difficult, but one has to be careful.  We address this
matter in \cite{W6}.

\vskip 1 cm

\begin{tabular}{l}
Department of Mathematics \\
University of California \\
Berkeley, California 94720--3840, U.S.A. \\
\\
{\tt jawolf@math.berkeley.edu}
\end{tabular}

\enddocument
\end